%% file: ep_sem_all_review.tex
\documentclass[nopreprintline]{elsarticle}

\input{packages.tex}
\input{macros.tex}

\begin{document}
    \journal{}

    \begin{frontmatter}

    \title{\papertitle}

    \author[mox]{P.~C.~Africa\corref{cor1}}
    \ead{pasqualeclaudio.africa@polimi.it}

    \author[mox]{M.~Salvador}
    \author[dicatam]{P.~Gervasio}
    \author[mox]{L.~Dede'}
    \author[mox,epfl]{A.~Quarteroni}

    \address[mox]{MOX -- Dipartimento di Matematica, Politecnico di Milano,\par P.zza Leonardo da Vinci, 32, 20133 Milano, Italy}
    \address[dicatam]{DICATAM, Università degli Studi di Brescia,\par Via Branze, 38, 25123 Brescia, Italy}
    \address[epfl]{Mathematics Institute, \'{E}cole Polytechnique F\'{e}d\'{e}rale de Lausanne,\par Av. Piccard, CH-1015 Lausanne, Switzerland (\textit{Professor Emeritus})}

    \cortext[cor1]{Corresponding author.}

    \begin{abstract}
    We propose a matrix--free solver for the numerical solution of the cardiac
    electrophysiology model consisting of the monodomain nonlinear
    reaction--diffusion equation coupled with a system of ordinary differential
    equations for the ionic species.
    Our numerical approximation is based on the
    high--order Spectral Element Method (SEM) to achieve accurate numerical
    discretization while employing a much
    smaller number of Degrees of Freedom than first--order Finite Elements.
    We combine vectorization with sum--factorization, thus allowing for a very
    efficient use of high--order polynomials in a high performance computing framework.
    We validate the effectiveness of our matrix--free solver in a variety of
    applications and perform different electrophysiological simulations ranging
    from a simple slab of cardiac tissue to a realistic four--chamber heart geometry.
    We compare SEM to SEM with Numerical Integration (SEM--NI), showing that they provide comparable results in terms of accuracy and efficiency.
    In both cases, increasing the local polynomial degree $p$ leads to better numerical results and smaller computational times than reducing the mesh size $h$.
    We also implement a matrix--free Geometric Multigrid preconditioner that
    results in a comparable number of linear solver iterations with respect to
    a state--of--the--art matrix--based Algebraic Multigrid preconditioner.
    As a matter of fact, the matrix--free solver proposed here yields up to
    45$\times$ speed--up with respect to a conventional matrix--based solver.
    \end{abstract}


    \end{frontmatter}

    \section{Introduction}

    Mathematical and numerical modeling of cardiac electrophysiology provides meaningful tools to address clinical problems \textit{in silico}, ranging from the cellular to the organ scale \cite{Gerach_2021, Gray_2018, Piersanti_2022, Potse_2014, Strocchi_2020}.
    For this reason, several mathematical models and methods have been designed to perform electrophysiological simulations \cite{Quarteroni_2019, Trayanova_2011}.
    Among these, we consider the monodomain equation coupled with suitable ionic models, which describes the space--time evolution of the transmembrane potential and the flow of chemical species across ion channels \cite{ColliFranzone_2014}.

    This set of combined partial and ordinary differential equations describes solutions that resemble those of a wavefront propagation problem, \textit{i.e.} manifesting very steep gradients.
    Despite being extensively used \cite{Arevalo_2016, Bayer_2022,
    Gillette_2021, MendoncaCosta_2019}, the Finite Element Method (FEM) with first order polynomials
    does not seem to be the most suitable to properly capture
    the physical processes underlying cardiac electrophysiology \cite{quarteroni_1994}.
    Indeed, in such cases, a very fine mesh resolution is required to obtain fully convergent numerical results \cite{Woodworth_2021}, which calls for an overwhelming computational burden.

    High--order numerical methods come into play to tackle this specific issue:
    Spectral Element Method (SEM) \cite{patera_1984,maday_1989,chqz07},
    high--order Discontinuous Galerkin (DG) \cite{Arnold-2001,Cockburn-1998}, Finite Volume Method (FVM)
    \cite{Leveque_2002}, or Isogeometric Analysis (IGA) \cite{chb_iga_book}
    account for small numerical dispersion and dissipation errors while allowing
    for converging numerical solutions with less Degrees of Freedom (DOFs)
    \cite{Bucelli_2021, Cantwell_2014, Coudiere_2017, Hoermann-2018}.
    However, the use of high--order polynomials in \textit{matrix--based} solvers for complex scenarios has been hampered by several numerical challenges, which are mostly related to the stiffness of the discretized monodomain problem \cite{Vincent_2015}.

    In this context, we develop and implement a high--order \textit{matrix--free} numerical solver
    that can be readily employed for CPU--based, massively
    parallel, large--scale numerical simulations.
    Since there is no need to assemble any matrix, all the floating point
    operations are associated with matrix--vector products that represent the most
    demanding computational kernels at each iteration of iterative solvers.
    Thanks to vectorization \cite{Arndt_2020}, which enables algebraic operations
    on multiple mesh cells at the same time, and sum--factorization
    \cite{Orszag_1980,Melenk_2001}, the higher the polynomial degree, the higher the computational advantages provided by the matrix--free solver \cite{Cantwell_2014,Kronbichler_2012}.
    Moreover, the small memory occupation required by the matrix--free implementation allows for its exploitation in GPU--based cardiac solvers \cite{Xia_2015,kronbichler_multigrid,DelCorso_2022}.

    In this manner, we obtain very accurate and efficient numerical simulations
    for cardiac electrophysiology, even if the linear solver remains unpreconditioned.
    Additionally, we implement a matrix--free Geometric Multigrid (GMG)
    preconditioner that is optimal for the values of $h$ (mesh size) and $p$
    (polynomial degree) considered in this paper when continuous model properties (\textit{i.e} a single ionic model and a continuous set of conductivity coefficients) are employed throughout the computational domain.

    We present different benchmark problems of increasing complexity for cardiac
    electrophysiology, ranging from the \textit{Niederer benchmark} on a slab of cardiac tissue
    \cite{niederer2011} to a whole--heart numerical simulation.
    We focus on two high--order discretization methods, namely, we compare SEM to
    SEM with Numerical Integration (SEM--NI), following the notations introduced in \cite{chqz07}.
    These two methods differ in the use of quadrature formulas, namely
    Legendre--Gauss for SEM and Legendre--Gauss--Lobatto for SEM--NI.
    Numerical results of Section~\ref{sec:numres} show that the two methods feature
    a similar behaviour in terms of both accuracy and computational costs. In both cases, choosing a higher polynomial degree $p$
    leads to a fairly more beneficial ratio between accuracy and computational
    costs than reducing the mesh size $h$. For instance,
    working with two discretizations with the same number of DOFs on
    the Niederer benchmark,
    the solution computed with ${\mathbb Q}_4$ (local polynomials of degree 4 with
    respect to each spatial coordinate) and average mesh size
    $h_{\mathrm{avg}}=0.48$ mm is more accurate than the one obtained with ${\mathbb
    Q}_1$ and average mesh size $h_{\mathrm{avg}}=0.12$ mm.
    Moreover, the former one has been computed at a computational cost that is about 40\% of the latter one.

    We also evaluate the performance of our matrix--free solver:
    a 45$\times$ speed--up is achieved with respect to the matrix--based solver.
    Furthermore, while with the matrix--based implementation the assembling and solving phases of the monodomain problem take more
    than 70\% of the total computational time, which also includes the solution of the system of ODEs associated with the coupled ionic models and the evaluation of the ionic current at each time step, plus some negligible initialization stages, this value
    drops to approximately 20\% with the matrix--free solver.

    The mathematical models and the numerical methods contained in this paper have been
    implemented in \lifex{} \cite{africa2022lifex} (\url{https://lifex.gitlab.io/}),
    a high-performance \texttt{C++} library developed within the iHEART project and based
    on the \texttt{deal.II} (\url{https://www.dealii.org}) Finite Element core \cite{dealII92}.

    The outline of the paper is as follows. We describe the monodomain model in Section~\ref{sec:problem}.
    We address its space and time discretizations in Section~\ref{sec:discretization}.
    We propose the matrix--free solver for cardiac electrophysiology and the matrix--free GMG preconditioner in Section~\ref{sec:mf-mb}, discussing details about vectorization, sum--factorization and highlighting similarities and differences between the matrix--based and the matrix--free solvers.
    Finally, the numerical results in Section~\ref{sec:numres} demonstrate the high
    efficiency of our high--order SEM matrix--free solver against the low--order
    FEM matrix--based one.

    \section{Mathematical model}\label{sec:problem}

    For the mathematical modeling of cardiac electrophysiology, we consider the
    monodomain equation coupled with suitable ionic models
    \cite{Quarteroni_2019,ColliFranzone_2014}:
    \begin{eqnarray}\label{monodomain}
    \left\{\begin{array}{ll}
    \displaystyle
    \Chi\left(\Cm\frac{\partial \Pot}{\partial t}+\Iion(\Pot,\Gating,\Conc)\right)
    -\nabla\cdot(\DiffTens\nabla \Pot)=\Chi\Iapp({\bf x},t),& \mbox{ in }\Omega\times(0,T],\\[2mm]
    (\DiffTens\nabla \Pot)\cdot {\bf n}=0, & \mbox{ on }\partial\Omega\times(0,T],\\[2mm]
    \displaystyle
    \frac{d\Gating}{dt}=\RhsGating(\Pot,\Gating,\Conc),
    & \mbox{ in }\Omega\times(0,T],\\[2mm]
    \displaystyle
    \frac{d\Conc}{dt}=\RhsConc(\Pot,\Gating,\Conc),
    & \mbox{ in }\Omega\times(0,T],\\[2mm]
    \Pot({\bf x},0)=\Pot_0({\bf x}),\
    \Gating({\bf x},0)=\Gating_0({\bf x}),\
    \Conc({\bf x},0)=\Conc_0({\bf x}), &  \mbox{ in }\Omega.
    \end{array}\right.
    \end{eqnarray}
    The unknowns are: the transmembrane potential $\Pot$, the vector
    $\Gating=(w_1,\ldots,w_M)$ of the probability density
    functions of $M$ gating variables, which represent the fraction of open channels across the membrane of a single cardiomyocyte, and the vector $\Conc
    =(z_1,\ldots, z_P)$ of the concentrations of $P$ ionic species.
    For the sake of simplifying the notation, in the following the membrane capacitance per unit area $\Cm$ and the membrane surface--to--volume ratio $\Chi$ are set equal to $1$.

    The mathematical expressions of the functions $\RhsGating(\Pot,\Gating,\Conc)$ and $\RhsConc(\Pot,\Gating,\Conc)$,
    which describe the dynamics of gating variables and ionic concentrations respectively,
    and the ionic current $\Iion(\Pot,\Gating,\Conc)$ strictly depend on the choice of the ionic model.
    Here, the TTP06 \cite{tentusscher_2006} ionic model is adopted for the slab and ventricular geometries, while the CRN \cite{courtemanche1998} ionic model is employed for the atria.
    The action potential is triggered by an external applied current $\Iapp({\bf x},t)$.

    The diffusion tensor $\DiffTens$ is expressed as follows

    \begin{equation}
    \DiffTens = \sigmal \fZero \otimes \fZero + \sigmat \sZero \otimes \sZero + \sigman \nZero \otimes \nZero,
    \label{eqn: diffusiontensor}
    \end{equation}
    where the vector fields $\fZero$, $\sZero$ and $\nZero$ express the fiber, the
    sheetlet and the sheet--normal (cross--fiber) directions, respectively \cite{Piersanti_2021,Africa_2022_fibers}.
    We also define longitudinal, transversal and normal conductivities as $\sigmal, \sigmat, \sigman \in \mathbb{R}^+$, respectively \cite{Piersanti_2021}. Homogeneous Neumann boundary conditions are prescribed on the whole boundary
    $\partial\Omega$ to impose the condition of electrically isolated domain, $\mathbf{n}$ being the outward unit normal vector to the boundary.

    In this paper, the computational domain $\Omega\subset {\mathbb R}^3$ is represented either by a slab of cardiac tissue or by the Zygote geometry \cite{Zygote}.

    \section{Space and time discretizations}\label{sec:discretization}

    In order to discretize in space the system (\ref{monodomain}), we adopt SEM \cite{patera_1984,maday_1989,chqz07,bernardi2001,karniadakis_sherwin},
    a high-order method that can be recast in the framework of the Galerkin method \cite{quarteroni_1994}.

    We consider a family of hexahedral conforming meshes,
    satisfying standard assumption of regularity and quasi--uniformity
    \cite{quarteroni_1994}, and let $h>0$ denote the mesh size.

    At each time, we look for the discrete
    solution belonging to the space of globally continuous functions
    that are the tensorial product of univariate piecewise (on each mesh element)
    polynomial functions of local degree $p\ge 1$ with
    respect to each coordinate. The local finite
    element space is referred to as ${\mathbb Q}_p$, while
    we denote by $V_{h,p}$ the global finite dimensional space.

    When using SEM, the univariate basis functions are of Lagrangian (\textit{i.e.},
    nodal) type and their
    support nodes $x_i$ are the Legendre--Gauss--Lobatto quadrature nodes
    (see, \textit{e.g.}, \cite[Ch. 2]{chqz06}), suitably mapped from the reference interval
    $[-1,1]$ to the local 1D elements.

    One of the main features of SEM is that, when the data are smooth enough,
    the induced approximation error decays more than algebraically fast with
    respect to the local polynomial degree. Indeed, it is said that SEM features exponential or spectral
    convergence.
    At the same time, the convergence with respect to the mesh size $h$ behaves
    as in FEM. More precisely, if $u\in H^s(\Omega)$, with $s > \frac{3}{2}$, denotes the exact solution of a linear second--order elliptic
    problem in a Lipschitz domain $\Omega$ and $u_\text{SEM}\in V_{hp}$ is its SEM approximation, the following error estimate holds

    \begin{equation}\label{sem_error}
    \|u-u_\text{SEM}\|_{H^1(\Omega)}\leq C  h^{\min(p+1,s)-1}p^{1-s}\|u\|_{H^s(\Omega)}.
    \end{equation}
    We refer, e.g., to \cite{chqz07,gervasio-2020} for an experimental support to this
    estimate.

    SEM can be considered as a special case of $hp-$FEM
    (\cite{karniadakis_sherwin,szabo_babuska_1991,schwab_1998}) with nodal basis
    functions and conforming hexahedral meshes.

    Typically, when using SEM, the integrals appearing in the Galerkin
    formulation of the differential problem \eqref{monodomain} are computed by the
    composite Legendre--Gauss (LG) quadrature
    formulas (see \cite{chqz07, chqz06}).  In principle, one can choose LG
    formulas of the desired order of exactness to guarantee a highly accurate
    computation of all the integrals appearing in \eqref{monodomain}.
    However, a typical choice is to use LG formulas with $(p+1)$ quadrature nodes,
    which guarantees that the entries of both the mass matrix and the stiffness
    matrix with constant coefficients are computed exactly while keeping the
    computational costs not too large \cite{Kronbichler_2012,Fehn_2019}.

    A considerable improvement in reducing the computational times of evaluating
    the integrals consists of using Legendre--Gauss--Lobatto (LGL) quadrature formulas (instead of LG ones),
    again with $(p+1)$ nodes that now coincide with the support nodes of the
    Lagrangian basis functions. This results into the so--called SEM--NI method (NI standing for Numerical
    Integration).
    Since the Lagrangian basis functions are mutually orthogonal with respect to the discrete
    $L^2-$inner product induced by the LGL formulas,
    the mass matrix of the SEM--NI method is diagonal, although not integrated exactly;
    this is a great strength of SEM--NI in solving time--dependent differential problems through explicit methods when the mass matrix is assembled.
    On the other hand, as the degree of exactness of LGL quadrature formulas
    using $(p+1)$ nodes along each direction is $2p-1$,
    the integrals associated with the nonlinear terms of the differential problem
    may introduce quadrature errors and aliasing effects that are as significant as
    the nonlinearities.

    We remark that ${\mathbb Q}_1-$SEM is equivalent to ${\mathbb Q}_1-$FEM,
    while ${\mathbb Q}_1-$SEM--NI is in fact ${\mathbb Q}_1-$FEM in which the
    integrals are approximated by the trapezoidal quadrature rule \cite{quarteroni_1994}.

    We choose the same local polynomial degree $p$ (and then the same
    finite dimensional space) for approximating the transmembrane potential $\Pot$, the
    gating variables $w_i$ (for $i=1,\ldots, M$) and the ionic concentrations $z_i$
    (for $i=1,\ldots, P$) at each time $t\in (0,T]$.

    All the time derivatives appearing in Equation~\eqref{monodomain}
    have been approximated using the 2nd--order Backward Differentiation Formula (BDF2)
    over a discrete set of time steps $t^n = n \Delta t, \,n = 0, \dots, N$, being $\Delta t$ the time step size.

One may wonder whether the BDF2 scheme is accurate enough for our simulations, even
when high-order spatial discretizations (${\mathbb Q}_3$ and ${\mathbb Q}_4$
SEM) are used, or if it is better to consider higher order methods like, \textit{e.g.}, the BDF3 scheme.
To remove any doubt, we have approximated the heat equation in a
two--dimensional domain, with discretization parameters similar to those used
in our simulations. We have verified that, with these discretization
parameters, the errors in space overbear those
in time, thus making the use of BDF3 worthless.
Moreover, we remark that, while BDF2 is absolutely stable, BDF3 is
not, then special care should be given to the choice of the time step. We refer
to \ref{Appendix} for a more in--depth analysis.

    Regardless of the quadrature formula (LG or LGL), the algebraic
    counterpart of the monodomain problem (\ref{monodomain}) reads:
    given $\AGating^0$, $\AConc^0$ and $\APot^0$, and suitable initializations for $\AGating^1$, $\AConc^1$ and $\APot^1$, then, for any $n\geq 1$, find $\AGating^{n+1}$, $\AConc^{n+1}$ and $\APot^{n+1}$ by solving the following partitioned scheme:

    \begin{subequations}
   	\label{electrophysiology_discr}
    \begin{align}
    \left\{\begin{aligned}
    & \frac{3\AGating^{n+1}-4\AGating^n+\AGating^{n-1}}{2\Delta t}=
    \RhsGating(\APot^{*},\AGating^{n+1},\AConc^{n+1}),\\
    &\frac{3\AConc^{n+1}-4\AConc^n+\AConc^{n-1}}{2\Delta t}=
    \RhsConc(\APot^{*},\AGating^{n+1},\AConc^{n+1}),
    \end{aligned}\right.\label{ionic_discr}\\[2mm]
    \Mass \frac{3\APot^{n+1}-4\APot^n+\APot^{n-1}}{2\Delta t}+ \Stiff
    \APot^{n+1} = \AIion^{n+1}+ \AIapp^{n+1}.\label{monodomain_discr}
    \end{align}
    \end{subequations}

    The arrays $\APot^{n+1}$ and $\AGating^{n+1}$ and $\AConc^{n+1}$
    contain the SEM or SEM--NI DOFs of the transmembrane potential, gating variables and ionic concentrations, respectively,
    $\Mass$ and $\Stiff$ are the SEM or SEM--NI mass and
    stiffness matrices, respectively, and $\AIapp_i^{n+1}=\Iapp({\bf x}_i,t^{n+1})$.
    The entries of $\AIion^{n+1}$ are computed with Ionic Current Interpolation
    (ICI) \cite{Regazzoni_2022}, \textit{i.e.},

    \begin{equation}\label{algebraic_ionic}
    \AIion_i^{n+1}= - \int_\Omega \left(\sum_j\Iion(\APot_j^{*},
    \AGating_j^{n+1}, \AConc_j^{n+1})\varphi_j\right)\varphi_i,
    \end{equation}
    with $\varphi_i$ the $i^{th}$ Lagrange basis function of the finite dimensional
    space $V_{h,p}$.  We remark that, when SEM--NI with LGL quadrature formulas are
    employed, ICI coincides with Lumped--ICI \cite{quarteroni2017},
    as the lumping of the SEM mass matrix coincides with the SEM--NI mass matrix.

    If we set $\APot^{*}=\APot^{n+1}$, then we recover the fully implicit BDF2 scheme. Nevertheless,
    we highlight that the function $\Iion$ is typically strongly nonlinear.
    To overcome the drawbacks of this nonlinearity, we adopt the extrapolation formula $\APot^{*}=2\APot^{n}-\APot^{n-1}$
    of  $\APot^{n+1}$, that is second--order accurate with respect to $\Delta t$.
    The resulting semi--implicit
    scheme is 2nd--order accurate in time when $\Delta t \to 0$ (see, \textit{e.g.},
    \cite{gervasio2006}).

    The  ordinary differential equations \eqref{ionic_discr} are
    associated with the ionic model and provide both the gating variables and the
    ionic species, while Equation~\eqref{monodomain_discr} is the
    discretization of the monodomain equation and its solution at the generic time
    step $n\geq 1$ is obtained by solving the linear system

    \begin{equation}\label{linsys}
    \Matrix\APot^{n+1}=\ARhs^{n+1},
    \end{equation}
    where

    \begin{equation}\label{matrix_rhs}
    \Matrix=\frac{3}{2\Delta t}\Mass +\Stiff, \hskip 1.cm
    \ARhs^{n+1}=\frac{1}{2\Delta
    t}\Mass(4\APot^n-\APot^{n-1})+\AIion^{n+1}+\AIapp^{n+1}.
    \end{equation}
    Solving the linear system (\ref{linsys})
    represents the most computationally demanding part of Equation~\eqref{electrophysiology_discr}.
    We refer to Section~\ref{sec:numres}, in particular to Table
    \ref{tab: matrix--free_matrix--based_percentages}, for further details about this specific aspect.

    \section{Matrix--free and matrix--based solvers}\label{sec:mf-mb}

    As in FEM, the matrix $\Matrix$ based on either SEM or SEM--NI has a very
    sparse structure, thus iterative methods are the natural candidates to solve the
    linear system (\ref{linsys}).
    Since $\Matrix$ is symmetric and positive definite, we have adopted
    the Conjugate Gradient (CG) method or its preconditioned version (PCG).

    Excluding the preconditioner step, the most expensive part of one CG--iteration
    is the evaluation of a matrix--vector product $\Matrix \mathbf{v}$, where
    $\mathbf{v}$ is a given vector.

    Typically, in a conventional matrix--based solver,
    the matrix $\Matrix$ is assembled and stored in sparse format, then referenced whenever the matrix--vector product has to be evaluated, \textit{i.e.} during each CG iteration.
    The matrix--based solver aims at minimizing the number of floating point
     operations required for such evaluation and is
    a winning strategy in ${\mathbb Q}_1-$FEM discretization for which the band of
    the matrix $\Matrix$ is small.

    When SEM or SEM--NI discretizations of local degree $p$ are employed,
    each cell counts $(p+1)^3$ DOFs. It follows that
    the typical bandwidth of SEM (or SEM--NI) stiffness matrices is about
    $C(p+1)^3$ (where $C$ is the maximum number of cells sharing one node of the
    mesh) and it exceeds widely
    that of ${\mathbb Q}_1-$FEM stiffness matrices.
    The large bandwidth of the SEM matrix $\Matrix$ can worsen the computational times of accessing the matrix entries, thus deteriorating the efficiency of the iterative solver.

    Moreover, in modern processors, access to the main memory has become the
bottleneck in many solvers for partial differential equations: a matrix--vector
product based on matrices requires far more time waiting for data to arrive
from memory than on actually doing the floating point operations. Thus, it is
demonstrated to be more efficient to recompute matrix entries -- or rather, the action of the differential operator represented by these entries on a known vector, cell by cell -- rather than looking up global matrix entries in the memory, even if the former approach requires a significant number of additional floating point operations \cite{Kronbichler_2012}.

    This approach is referred to as \emph{matrix--free}. In practice, shape
    functions values and gradients are pre-computed for each basis function on the reference cell, for each quadrature node. Then, the Jacobian of the transformation from the real to the reference cell is cached, thus improving the computational cost of the evaluation.

	\subsection{Vectorization and sum--factorization}
\begin{figure}
    	\centering
    	\includegraphics[width=0.3\textwidth]{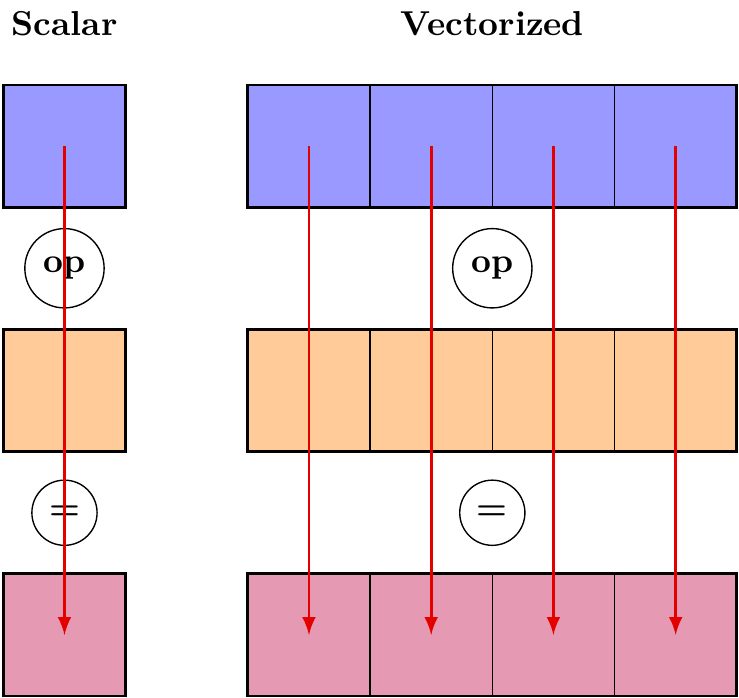}
    	\caption{Comparison between scalar and vectorized operations.}
    	\label{fig:vectorization}
    \end{figure}

\begin{figure}
        \centering
        \begin{subfigure}[t]{0.3\linewidth}
            \centering
            \includegraphics[width=\textwidth]{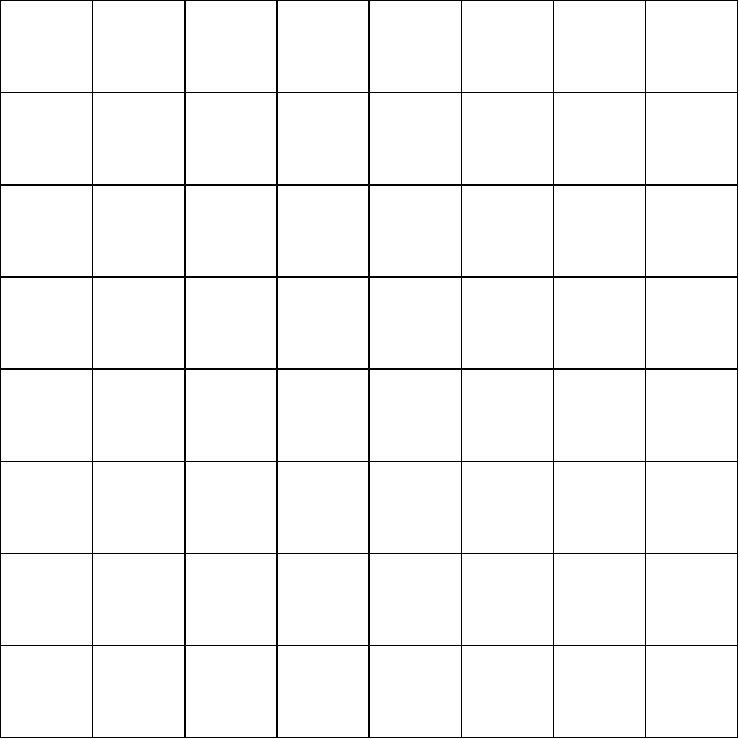}
            \caption{Original mesh.}
        \end{subfigure}
        \hfill
        \begin{subfigure}[t]{0.3\linewidth}
            \centering
            \includegraphics[width=\textwidth]{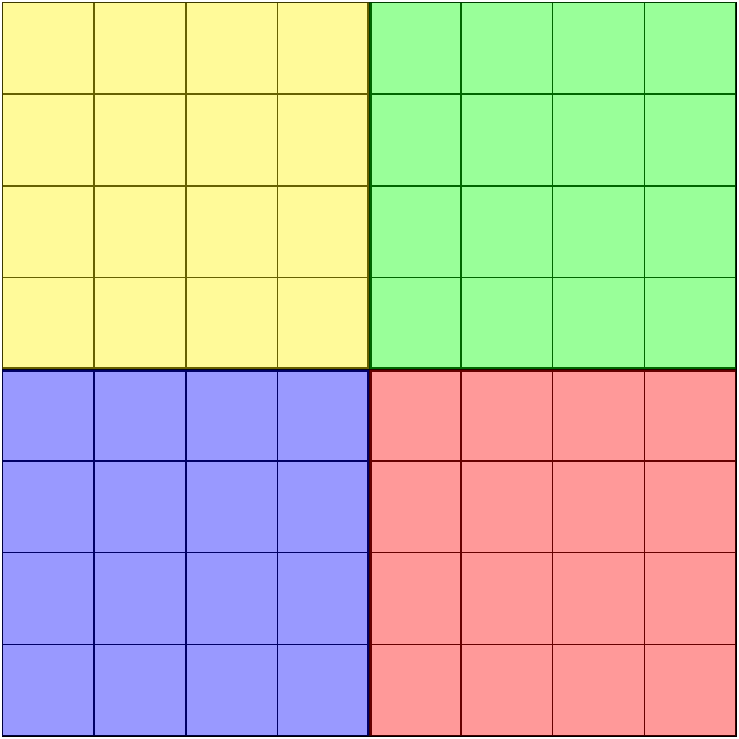}
            \caption{Parallel partitioning of mesh cells. Bold solid lines delimit cells owned by different parallel processors.}
        \end{subfigure}
        \hfill
        \begin{subfigure}[t]{0.3\linewidth}
            \centering
            \includegraphics[width=\textwidth]{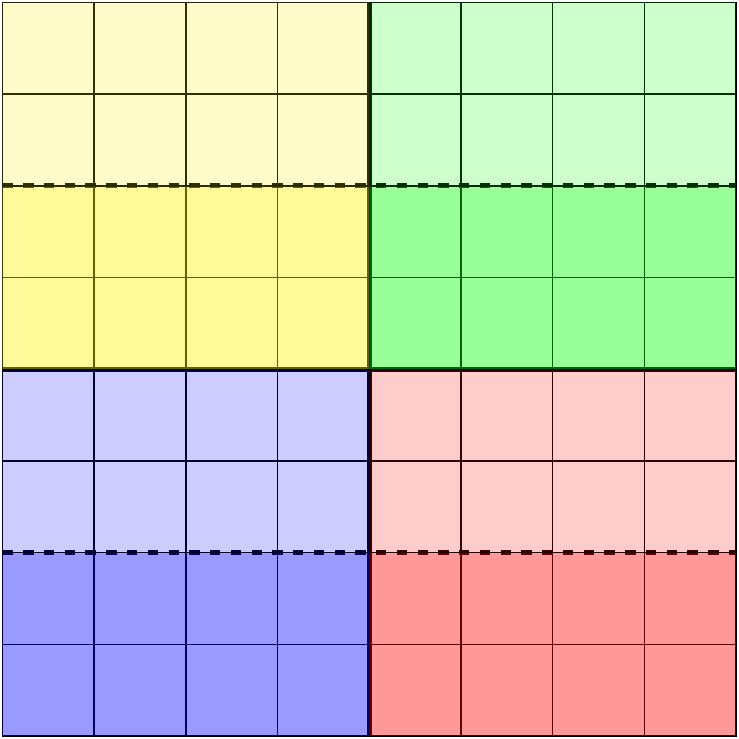}
            \caption{Vectorization of cell operations on each parallel unit. Dashed lines delimit batches of cells processed using vectorized operations, within each processor.}
        \end{subfigure}
        \caption{Vectorization in a parallel framework: when the original mesh is partitioned and distributed among multiple computational units, vectorization is applied at the level of each processor (in this example, 4 parallel processors are considered). Here, different colors refer to different parallel units and light/dark variations represent different vectorized batches (in this example, vectorization acts on 8 cells).}
        \label{fig:vectorization_mesh}
    \end{figure}

    In FEM solvers (and, similarly, in SEM ones), the cell--wise computations are
    typically exactly the same for all cells, and hence a Single--Instruction,
    Multiple--Data (SIMD) stream can be used to process several values at once (see Figure~\ref{fig:vectorization}). Vectorization is a SIMD concept, that is, one CPU instruction is used to process multiple cells at once. Modern CPUs support SIMD instruction sets to different extents, \textit{i.e.} one single CPU instruction can simultaneously process from two doubles (or four floats) up to eight doubles (or sixteen floats), depending on the underlying architecture \cite{cebrian2020scalability}. Additionally, vectorization can also be combined to a distributed memory parallelism \cite{zhong2022using}. In our case, the mesh cells and degrees of freedom are partitioned and distributed among different parallel processing units via MPI \cite{mpi40}, resulting in the scheme shown in Figure~\ref{fig:vectorization_mesh}.

	Vectorization is beneficial only in arithmetic intensive operations, whereas additional computational power becomes useless when the workload bottleneck is the memory bandwidth.
    For this reason, vectorization is typically not used explicitly in
\textit{matrix--based} Finite Element codes, whose computational efficiency is
dominated by memory access. On the other hand, \textit{matrix--free} solvers
can easily benefit from the additional computational speed--up brought in by
vectorized operations. As a matter of fact, in this case a matrix--vector product results from
recomputing the local action of the matrix on the vector, cell by cell, every
time that it is needed,
    rather than accessing global matrix entries in the memory.
    In our matrix--free Algorithm \ref{algo:matrix_based}, vectorization acts on the cell loop performed at line \ref{algo:vectorized_loop}.

	Finally, thanks to the fact that the multivariate SEM Lagrange basis is of
    tensorial type, in order to reduce the computational complexity of one
    evaluation of the product $\Matrix \mathbf{v}$, \emph{sum--factorization}
    can also be exploited \cite{Orszag_1980,Melenk_2001,Kronbichler_2019}. In
    this way, the matrix–vector product for the Laplace operator in the generic three dimensional cell
    requires only $9(p+1)$ floating point operations instead of $(p+1)^3$ per
    degree of freedom, resulting in a complexity equal to $\mathcal{O}\left(9(p+1)^4\right)$ instead of $\mathcal{O}\left((p+1)^6\right)$, and this
    still plays in favor of repeating computations rather than accessing the
    memory (see \cite[Section 2.3.1]{Cantwell_2011} and \cite[Section 4.5.1]{chqz06}).

    On these bases, a high--order matrix--free solver is more efficient both in
    terms of memory occupation (no system matrix is assembled and stored globally)
    and computational time \cite{Kronbichler_2012}, as we will also show in Section~\ref{sec:numres}.

	\subsection{Application to cardiac electrophysiology}
	\label{sec:mf_vs_mb}
	In Algorithms \ref{algo:matrix_free} and \ref{algo:matrix_based} we
	display the computational workflow resulting from applying the matrix--free and the matrix--based solver to problem \eqref{electrophysiology_discr}, respectively. The different phases are listed. We highlight that the operations needed to solve the ionic model \eqref{ionic_discr} and to compute \(\Iion\) are the same for both algorithms.

	In the following, the expression \textit{assembly phase} will refer to
    the assembly of the right--hand side, in the case of the matrix--free solver, and to the assembly of both the right--hand side and the system matrix, in the case of the matrix--based solver. On the other hand, the \textit{linear solver} phase of the matrix--free algorithm encloses also a cell loop for computing the local action of the discretized operator at each CG iteration, which is not required in the matrix--based case.

	Finally, we remark that the diffusion tensor \(\DiffTens\) is evaluated as in Equation~\eqref{eqn: diffusiontensor} at every quadrature point, resulting in $9(p+1)^3$ additional memory accesses per cell.

	Both Algorithms \ref{algo:matrix_free} and \ref{algo:matrix_based} are run in parallel as described in the previous section, by distributing all loops over DOFs and cells according to the mesh partitioning. For the sake of simplicity, the application of suitable preconditioners is omitted from the listings. More details on this topic are discussed in the next section.

	\begin{algorithm}
	\caption{Workflow for the matrix--free solver.}
	\label{algo:matrix_free}
	\KwIn{$\AGating^0$, $\AConc^0$ and $\APot^0$, and suitable initializations for $\AGating^1$, $\AConc^1$ and $\APot^1$}
	\ForEach{time step \(n \geq 1\)}{
		Compute \(\APot^{*}=2\APot^{n}-\APot^{n-1}\)\;
		\nonl\;
		\nonl Solve the ionic model and compute the ionic current DOF--wise:\;
		\ForEach{DOF \(i\)}{
			Solve Eq.~\eqref{ionic_discr} for \(\AGating_i^{n+1}\) and \(\AConc_i^{n+1}\)\;
			Compute \(\Iion_{,i} = \Iion(\APot_i^{*},
			\AGating_i^{n+1}, \AConc_i^{n+1})\)\;
		}
		\nonl\;
		\nonl \textbf{Assembly phase}: assemble the rhs of Eq.~\eqref{monodomain_discr}:\;
		\ForEach{mesh cell \(K\)}{
			Compute the local right--hand side\;
			\nonl \quad \(\ARhs_K = \AIion_K^{n+1} + \AIapp_K^{n+1} + \frac{\Mass_K}{2\Delta t}\left(4\APot_K^n - \APot_K^{n-1}\right)\)\;

			Compress \(\ARhs_K\) into a global right--hand side \(\ARhs\)\;
		}
		\nonl\;
		\nonl Solve the monodomain equation \(\left(\frac{3}{2\Delta t}\Mass + K\right)\APot^{(n+1)} = \ARhs\)\;
		\nonl \quad as in Eq.~\eqref{monodomain_discr} using CG:\;
		\ForEach{CG iteration \(m\) \textbf{until} convergence}{
			Determine the new conjugate vector \(\mathbf{v}\)\;
			\nonl\;
			\ForEach{mesh cell \(K\) \tcp*[r]{Vectorized loop.}\label{algo:vectorized_loop}}{
				Compute the local mass matrix \(\Mass_K\)\;
				Evaluate the diffusion tensor \(\DiffTens\) at quadrature nodes over \(K\)\;
				Compute the local stiffness matrix \(\Stiff_K\)\;
				Compute the local matrix--vector product:\;
				\nonl \quad \(\mathbf{z}_K = \left(\frac{3}{2\Delta t}\Mass_K + \Stiff_K\right)\mathbf{v}_K\)\;
				Compress the local contribution \(\mathbf{z}_K\) into a global vector \(\mathbf{z}\)\;
			}
			\nonl\;
			Apply the CG steps to \(\mathbf{z}\) in order to get the new solution \(\APot_m^{(n+1)}\)\;
		}
		\nonl\;
		\(\APot^{(n+1)} \gets \APot_m^{(n+1)}\)\;
	}
	\end{algorithm}

	\begin{algorithm}
		\caption{Workflow for the matrix--based solver.}
		\label{algo:matrix_based}
		\KwIn{$\AGating^0$, $\AConc^0$ and $\APot^0$, and suitable initializations for $\AGating^1$, $\AConc^1$ and $\APot^1$}
		\ForEach{time step \(n \geq 1\)}{
			Compute \(\APot^{*}=2\APot^{n}-\APot^{n-1}\)\;
			\nonl\;
			\nonl Solve the ionic model and compute the ionic current DOF--wise:\;
			\ForEach{DOF \(i\)}{
				Solve Eq.~\eqref{ionic_discr} for \(\AGating_i^{n+1}\) and \(\AConc_i^{n+1}\)\;
				Compute \(\Iion_{,i} = \Iion(\APot_i^{*},
				\AGating_i^{n+1}, \AConc_i^{n+1})\)\;
			}
			\nonl\;
			\nonl \textbf{Assembly phase}: assemble the global system matrix\;
			\nonl \quad and right--hand side:\;
			\ForEach{mesh cell \(K\)}{
				Compute the local mass matrix \(\Mass_K\)\;
				Evaluate the diffusion tensor \(\DiffTens\) at quadrature nodes over \(K\)\;
				Compute the local stiffness matrix \(\Stiff_K\)\;
				\nonl\;
				Compress the local system matrix \(\Matrix_K = \left(\frac{3}{2\Delta t}\Mass_K + \Stiff_K\right)\)\;
				\nonl \quad into a global matrix \(\Matrix\)\;
				Compress the local right--hand side\;
				\nonl \quad \(\ARhs_K = \AIion_K^{n+1} + \AIapp_K^{n+1} + \frac{\Mass_K}{2\Delta t}\left(4\APot_K^n - \APot_K^{n-1}\right)\)\;
				\nonl \quad into a global right--hand side \(\ARhs\)\;
			}
			\nonl\;
			Solve the monodomain equation \(\Matrix\APot^{(n+1)} = \ARhs\)\;
			\nonl \quad as in Eq.~\eqref{monodomain_discr} using CG\;
		}
	\end{algorithm}

	\subsection{A Geometric Multigrid matrix--free preconditioner}
	In order to precondition the CG method we have chosen Multigrid preconditioners.
	For the matrix--based solver, the Algebraic Multigrid (AMG) preconditioner \cite{janssen2011,Xu_2017} turns out to be a very efficient choice.
	Nevertheless, its implementation requires the explicit knowledge of the entries of the matrix $\Matrix$.

	Hybrid multigrid algorithms with matrix--free implementation
	for high--order discretizations have recently been proposed and discussed in \cite{Bastian-2019,Fehn-2020}. In particular the methods proposed in \cite{Fehn-2020} combine $h-$coarsening, $p-$coarsening, and AMG on the coarsest level, and they fully exploit the advantages of matrix--free algorithms with sum--factorization for the multigrid smoothers. The matrix assembly at the coarsest level is however required to implement the AMG solver.  We refer to \cite{Fehn-2020} for an interesting presentation of hybrid multigrid techniques and of the challenges to face for improving the efficiency of these algorithms.

	To overcome the drawback of assembling the matrix $\Matrix$ even at the
	coarsest level, we have adopted a fully Geometric
	Multigrid (GMG) preconditioner, more precisely the \textit{high--order $h$--multigrid} preconditioner \cite{Sundar-2015}, which uses $p$--degree interpolation and restriction among geometrically coarsened meshes. GMG methods are among the most efficient solvers for linear systems arising from the discretization of elliptic partial differential equations, offering an optimal complexity $\mathcal{O}(n)$ in the number of unknowns $n$, and they are often used as very efficient preconditioners (see \cite{kronbichler_multigrid, janssen2011, Trottenberg, Clevenger2021} and the literature cited therein).

\begin{figure}
		\centering
		\includegraphics[width=0.3\textwidth]{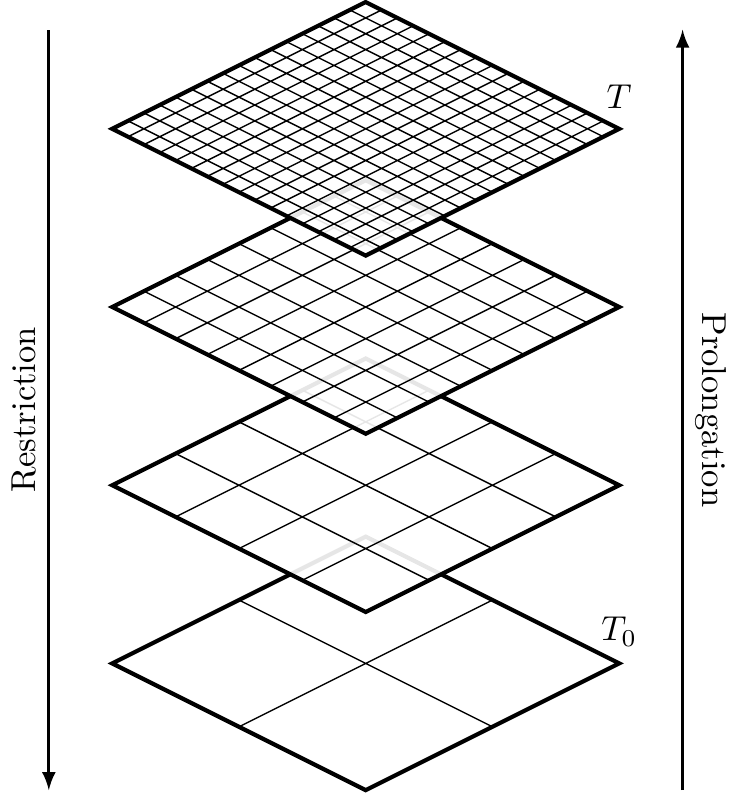}
		\caption{Schematization of multigrid methods in a two dimensional case. Starting from the real mesh \(T\), the action of the \textit{restriction} and \textit{prolongation} operators is shown, down to the coarse mesh \(T_0\).}
		\label{fig:multigrid}
	\end{figure}

	In the spirit of \cite{Kronbichler_2012,kronbichler_multigrid}, our GMG implementation relies
	on the simple yet effective scheme that considers a polynomial variant
	of the point--Jacobi smoother, namely a Chebyshev method with optimal parameters determined
	by an eigenvalue estimation based on Lanczos iterations \cite{Adams_2003}. This choice turns out
	to be very efficient in a matrix--free context because all its
  	computational kernels, including the smoother and the transfer
  	between different grid levels, are based on matrix--vector products involving
  	suitable collections of mesh cells \cite{Adams_2003}. In our case, a hierarchical collection of octree meshes is built by the recursive subdivision of each cell into 8 subcells, starting from a coarse mesh $T_0$ of size
  	$h_0$, as shown in Figure~\ref{fig:multigrid}.
  	Despite the higher throughput provided by multigrid in single precision \cite{Kronbichler_2012}, due to the high accuracy required by the monodomain problem at hand we decided to evaluate our GMG preconditioner in double precision, consistently with the matrix and vector representations within the CG solver.

    \section{Numerical results}\label{sec:numres}

    We present several numerical simulations of cardiac electrophysiology.
    First, we consider a benchmark problem on a slab of cardiac tissue \cite{niederer2011}, in order
    to compare SEM against SEM--NI and matrix--free against matrix--based in terms of computational efficiency and numerical accuracy.
    Then, we employ the Zygote left ventricle geometry \cite{Zygote} and we analyze the sole
    impact of increasing $p$, \textit{i.e.} the local polynomial degree, on the numerical solution.
    Finally, for the sake of completeness, we show the capability of our
    matrix--free solver by presenting a detailed Zygote four--chamber heart \cite{Zygote} electrophysiological simulation in sinus rhythm.

    For the time discretization, we use the BDF2 scheme with a time step $\Delta
    t=\SI{0.1}{\milli\second}$. The final time $T$ differs with the specific test
    case. We employ $T=\SI{0.2}{\second}$ in the Niederer benchmark \cite{niederer2011}, while $T=\SI{0.6}{\second}$ and $T=\SI{0.8}{\second}$ are considered for the left ventricle and whole--heart geometries, respectively.

	For what concerns the GMG--preconditioned CG solver in the matrix--free setting, we estimate the largest eigenvalue \(\tilde{\lambda}_\mathrm{max}\) for the Chebyshev smoother on each level by performing 10 CG sub--iterations, set the smoothing range to \([0.08\tilde{\lambda}_\mathrm{max},1.2\tilde{\lambda}_\mathrm{max}]\) (the number 1.2 is
	a safety factor that allows for some inaccuracies in the eigenvalue estimate) and choose a polynomial degree of 5 (\textit{i.e.}, 5 matrix–-vector products per level and iteration). Besides, for the AMG--preconditioned matrix--based CG solver we rely on the \texttt{Trilinos ML} smoothed aggregation \cite{trilinos_ml}, by performing 3 \(V\)--cycles with polynomial Chebyshev smoother of order \(3\) and by setting the aggregation threshold to \(10^{-1}\). All the parameters reported above have been empirically determined in order to keep the number of PCG iterations as low as possible.

	In all cases the PCG solver is run with a stopping criterion based on the absolute residual with tolerance $10^{-15}$.

	To compute the solution at time $t^{n+1}$, we use the solution at
	time $t^n$ as initial guess for the PCG
	algorithm.
	Because in our simulations we take a very small $\Delta t$ (tipically $\Delta
	t=10^{-4}$), the initial guess is
	itself a good approximation  of the solution and a low number of iterations is
	needed to satisfy the stopping criterion. We have verified that the norm of the
	starting residual of the linear system is between $10^{-9}$ and $10^{-8}$ along the whole numerical simulation, thus,
	an absolute stopping test on the residual with tolerance $10^{-15}$ means that
	we reduce the norm of the residual of about six to seven orders of magnitude.

    In the two test cases involving
    the slab and left ventricle, we employ the GMG (AMG) preconditioner for
    the matrix--free (matrix--based) solver.
    On the other hand, no preconditioner is
    introduced in the four--chamber heart numerical simulation, as the presence
    of different ionic models in the computational domain, namely the CRN model (\cite{courtemanche1998}) for atria and the TTP06 one (\cite{tentusscher_2006}) for ventricles, would make our GMG preconditioner non--optimal in \(h\) and \(p\).

    In Table~\ref{tab:parameters electrophysiology} we report the parameters of the monodomain equation. In particular, the conductivity tensor depends on the fiber distribution as in Equation~\eqref{eqn: diffusiontensor}, which is generated by means of the Laplace--Dirichlet Rule--Based
    Methods proposed in \cite{Piersanti_2021,Africa_2022_fibers}.

    \begin{table}[t!]
        \centering
        \begin{tabular}{l r r | l r r}
            \toprule
            Variable & Value & Unit & Variable & Value & Unit \\
            \midrule
            \multicolumn{3}{l|}{\textbf{Conductivity tensor}} & \multicolumn{3}{l}{\textbf{Applied current}} \\
            $\sigmal$ & \num{0.7643e-4} & \si{\meter\squared\per\second} &
            ${\widetilde{\mathcal{I}}_{\mathrm{app}}^{\mathrm{max}}}$  & \num{15} & \si{\volt\per\second}   \\
            $\sigmat$ & \num{0.3494e-4} & \si{\meter\squared\per\second} &
            $t_{\mathrm{app}}$  & \num{3e-3} & \si{\second}   \\
            $\sigman$ & \num{0.1125e-4} & \si{\meter\squared\per\second} \\
            \bottomrule
        \end{tabular}
        \caption{Parameters of the electrophysiology model. The conductivity tensor is defined as in Equation~\eqref{eqn: diffusiontensor}. For the CRN and
    TTP06 ionic models, we adopt the parameters reported in the original papers \cite{tentusscher_2006,courtemanche1998} for epicardial cells.}
        \label{tab:parameters electrophysiology}
    \end{table}

    The external current $\Iapp({\bf x},t)=\widetilde{\mathcal{I}}_{\mathrm{app}}^{\mathrm{max}}$ is
    applied for  $t \in (0, t_{\mathrm{app}}]$ in a cuboid for the Niederer benchmark
    (as described in \cite{niederer2011}), otherwise in different spheres for the ventricle and whole--heart
    test cases (we can deduce them from the numerical results shown in Figures \ref{fig:slab}, \ref{fig:activation_time_LV} and \ref{fig:solution_4CH}).

    All the numerical simulations were performed by using one cluster node endowed with 56 cores (two Intel Xeon Gold 6238R, 2.20 GHz), which is available at MOX, Dipartimento di Matematica, Politecnico di Milano.

    \subsection{Slab of cardiac tissue}
    \label{sec:slab}

    \begin{figure}
        \centering
        \begin{subfigure}{0.45\linewidth}
              \centering
              \includegraphics[width=0.85\textwidth]{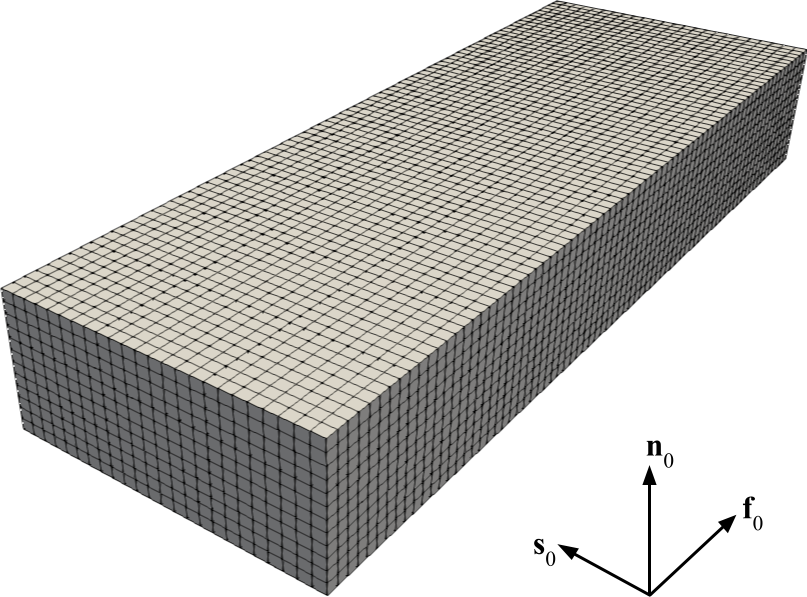}
        \end{subfigure}
        \begin{subfigure}{0.45\linewidth}
              \centering
              \includegraphics[width=\textwidth]{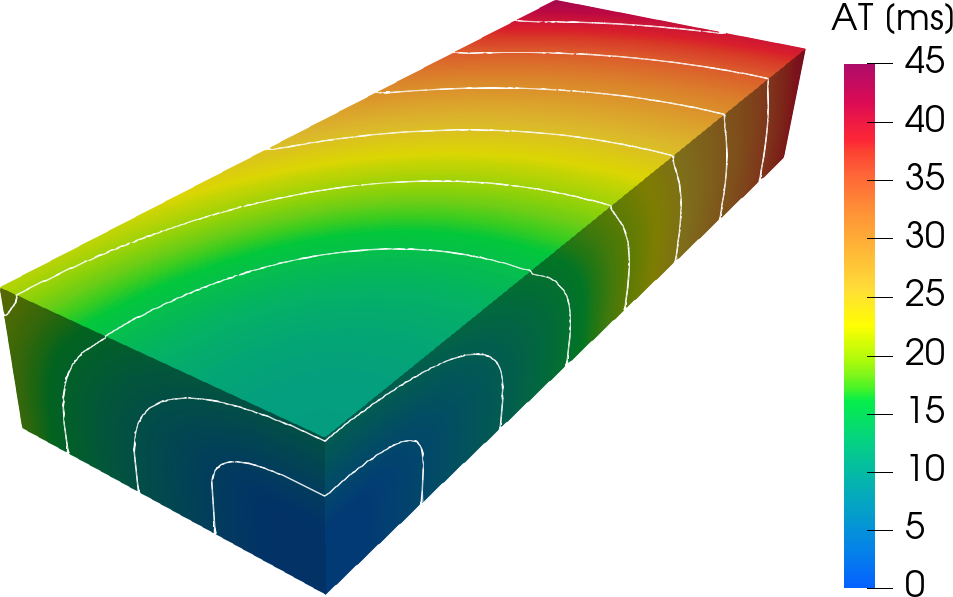}
        \end{subfigure}
        \caption{\emph{Niederer benchmark}. The geometry with the corresponding fiber, sheetlet and cross--fiber orientations and an example of mesh (left) and
    the associated simulated activation map (right). The blue color represents the region where the cubic stimulus is initially applied; the red one is associated with the corresponding diagonally opposite vertex, which is activated as last.}
        \label{fig:slab}
    \end{figure}

    The computational domain with an example of
    mesh (left) and the associated numerical simulation (right) for the Niederer
    benchmark \cite{niederer2011} is depicted in Figure~\ref{fig:slab}.
    An external stimulus of cubic shape is applied at one vertex, the electric signal propagates through the slab, and the diagonally opposite vertex is activated as the last point. The domain is discretized by means of a structured, uniform hexahedral mesh.

    We present a systematic comparison between SEM and SEM--NI for several values of both the mesh size $h$ and the local polynomial degree $p$, in order to understand which is the best formulation in terms of accuracy and computational cost.
    Moreover, we compare the efficiency of the matrix--free and matrix--based
    solvers for SEM.

\begin{figure}
      \centering
            \begin{subfigure}{.9\linewidth}
          \centering
          \includegraphics[keepaspectratio,
    width=\textwidth]{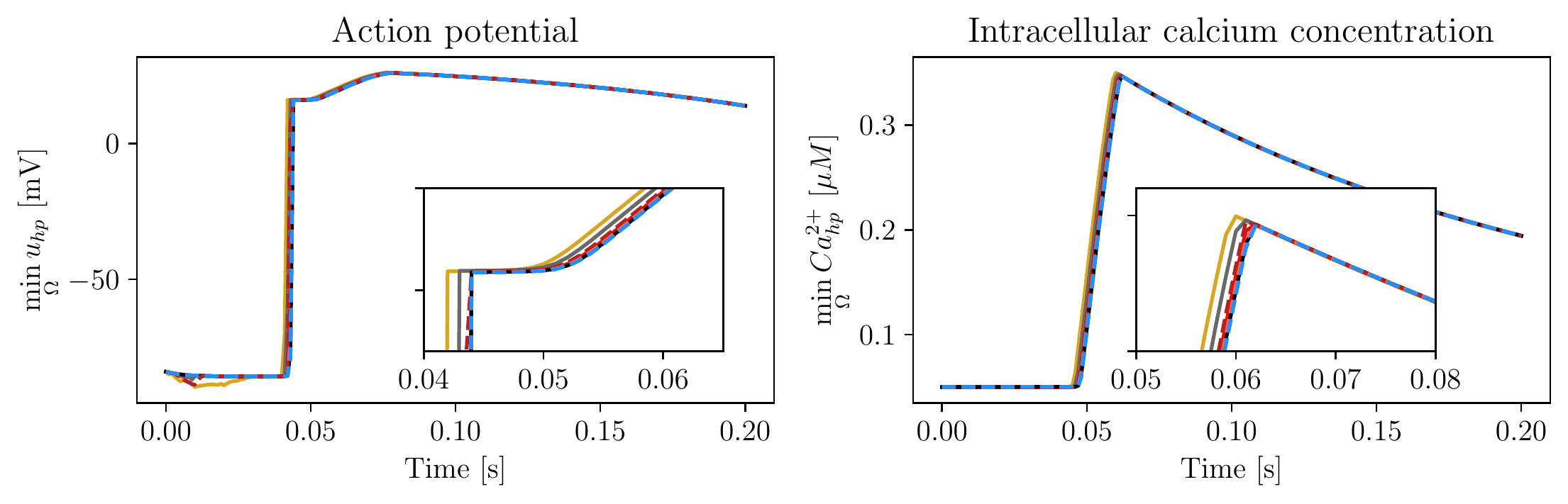}
        \end{subfigure}

        \hfill

        \begin{subfigure}{.9\linewidth}
          \centering
          \includegraphics[keepaspectratio,
    width=\textwidth]{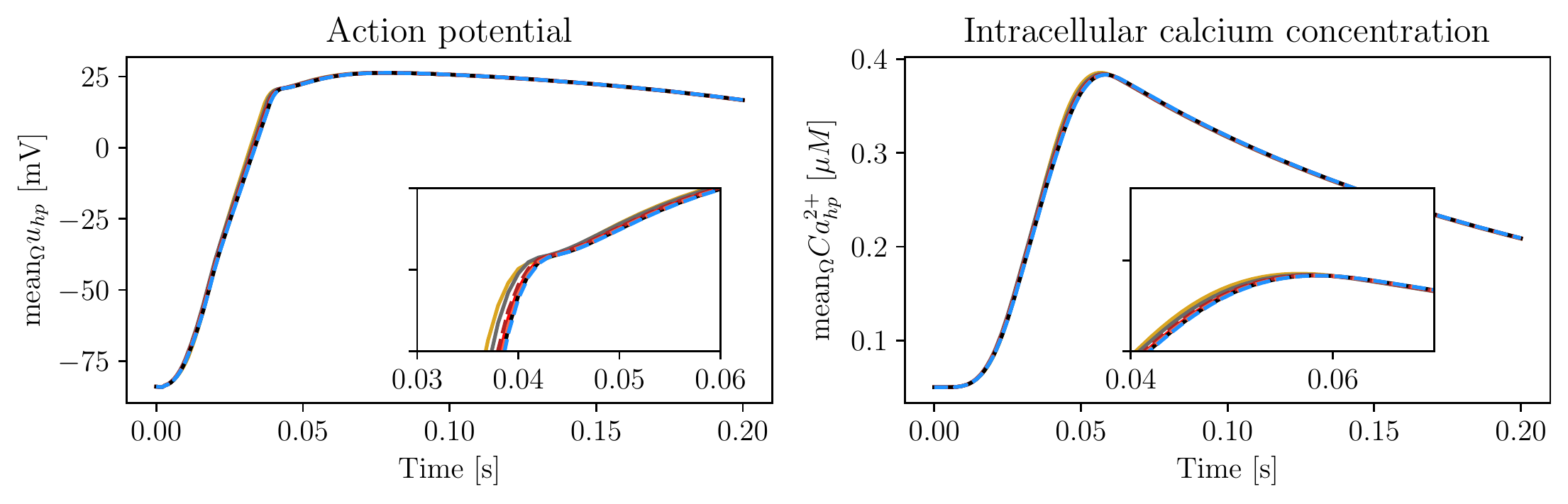}
        \end{subfigure}

        \hfill

            \begin{subfigure}{.9\linewidth}
          \centering
          \includegraphics[keepaspectratio,
    width=\textwidth]{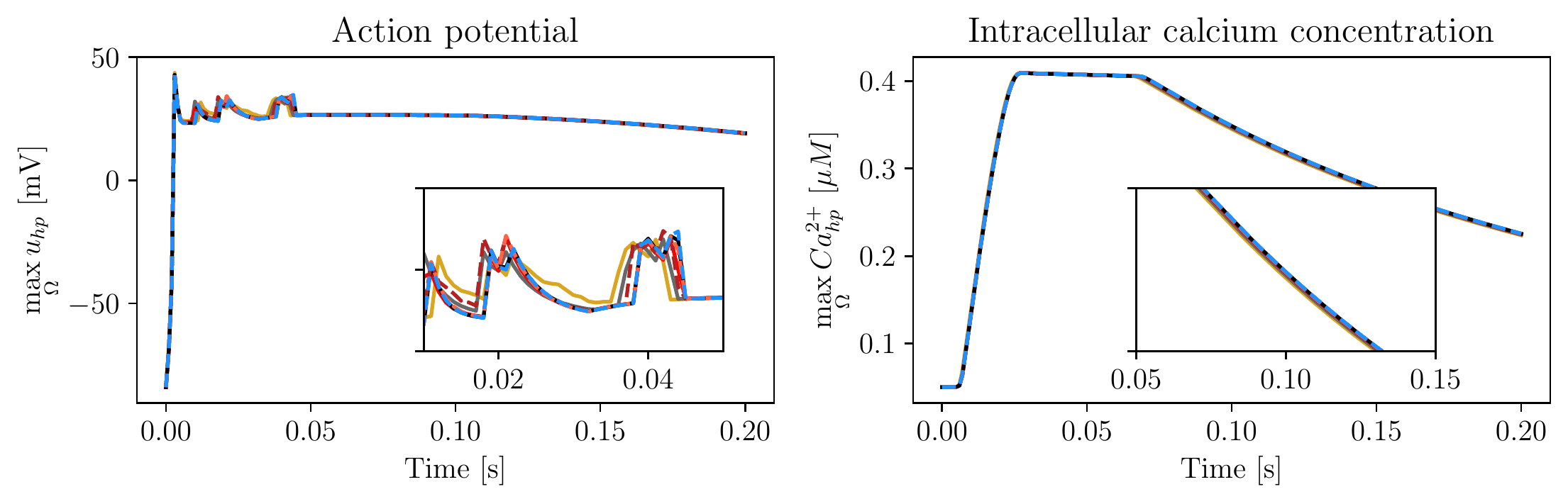}
       \end{subfigure}

        \hfill

            \begin{subfigure}{.9\linewidth}
          \centering
          \includegraphics[keepaspectratio,
    width=\textwidth]{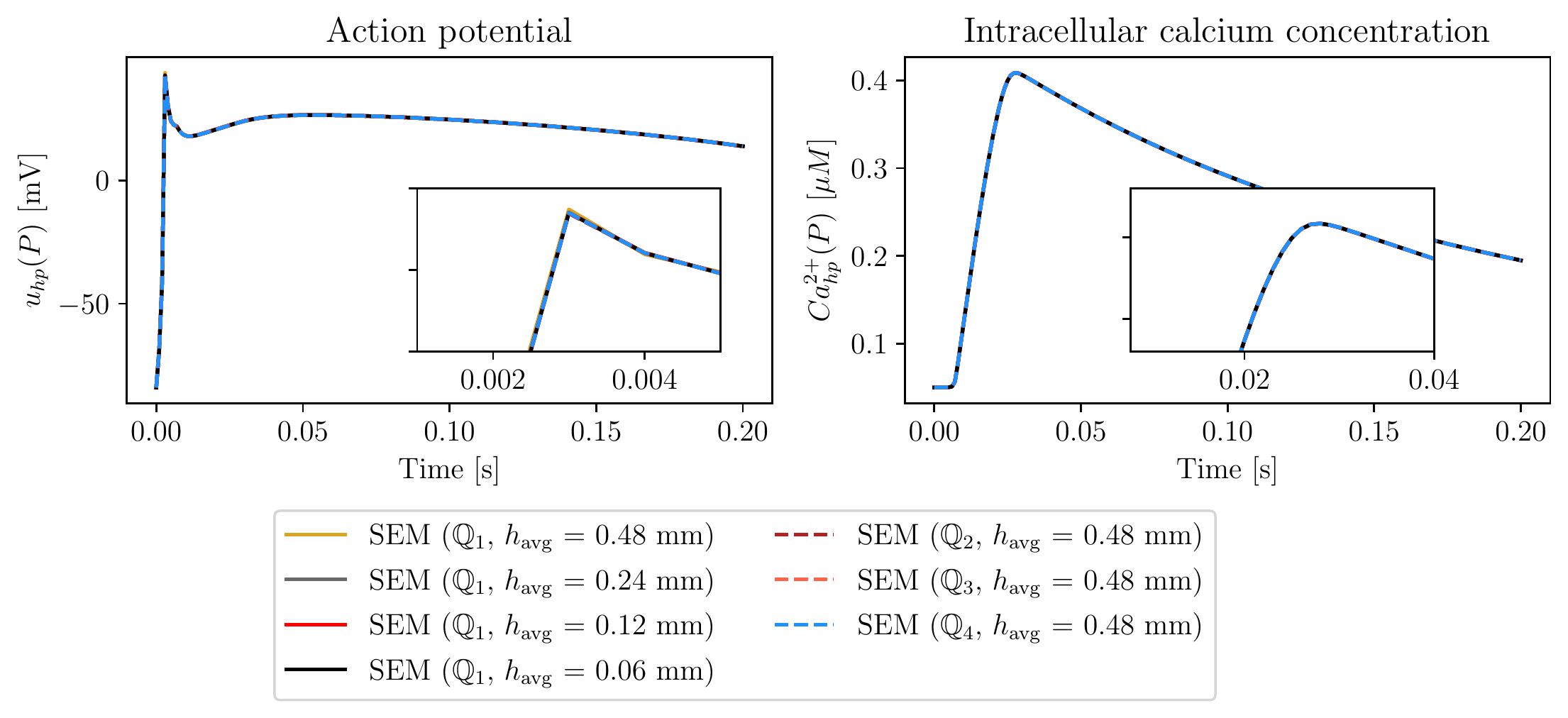}
       \end{subfigure}
      \caption{\emph{Niederer benchmark}. Minimum (top),
    average (second), maximum
    (third)  and point values (bottom) of the action
    potential $u$ and intracellular calcium concentration $Ca^{2+}$ over time for a
    slab of cardiac tissue. $P$ is a random point within the computational domain
    away from the initial stimulus. We consider different $hp$ combinations:
    SEM ${\mathbb Q}_1$ to ${\mathbb Q}_4$ and
    $h_\mathrm{avg} = \SI{0.48}{\milli\meter}$ to $h_\mathrm{avg} =
    \SI{0.06}{\milli\meter}$ ($h_\mathrm{avg}$ is the average mesh size).}
      \label{fig: slab_ap_calcium_sem}
    \end{figure}

\begin{figure}
      \centering
            \begin{subfigure}{.9\linewidth}
          \centering
          \includegraphics[keepaspectratio,
    width=\textwidth]{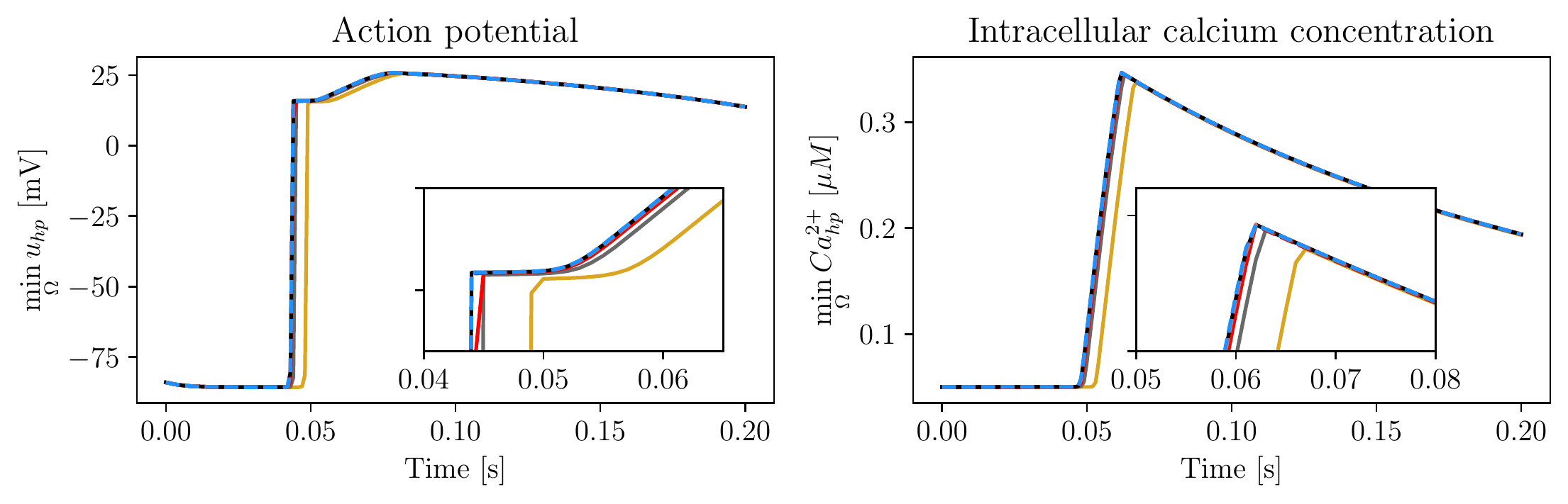}
        \end{subfigure}

        \hfill

        \begin{subfigure}{.9\linewidth}
          \centering
          \includegraphics[keepaspectratio,
    width=\textwidth]{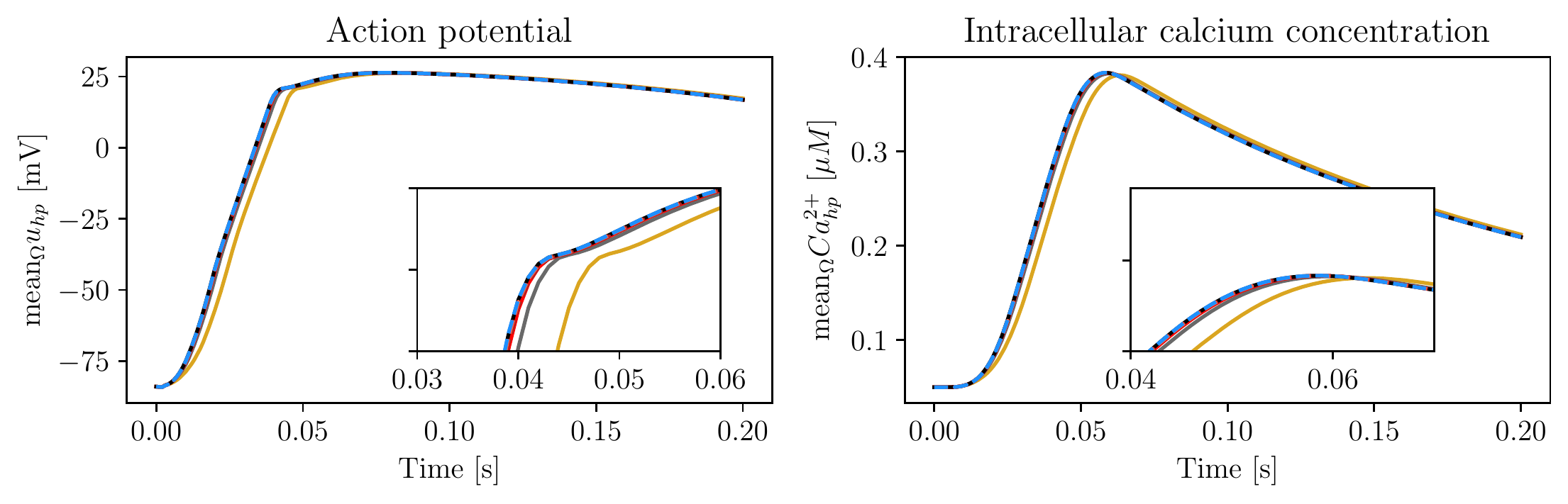}
        \end{subfigure}

        \hfill

            \begin{subfigure}{.9\linewidth}
          \centering
          \includegraphics[keepaspectratio,
    width=\textwidth]{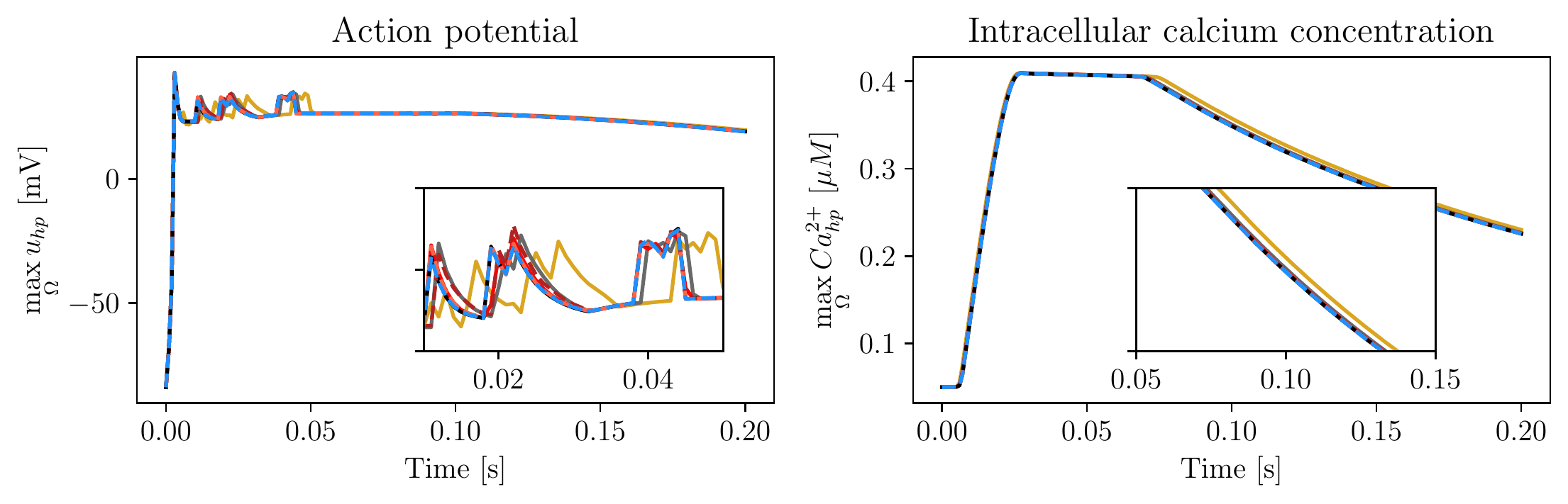}
       \end{subfigure}

        \hfill

            \begin{subfigure}{.9\linewidth}
          \centering
          \includegraphics[keepaspectratio,
    width=\textwidth]{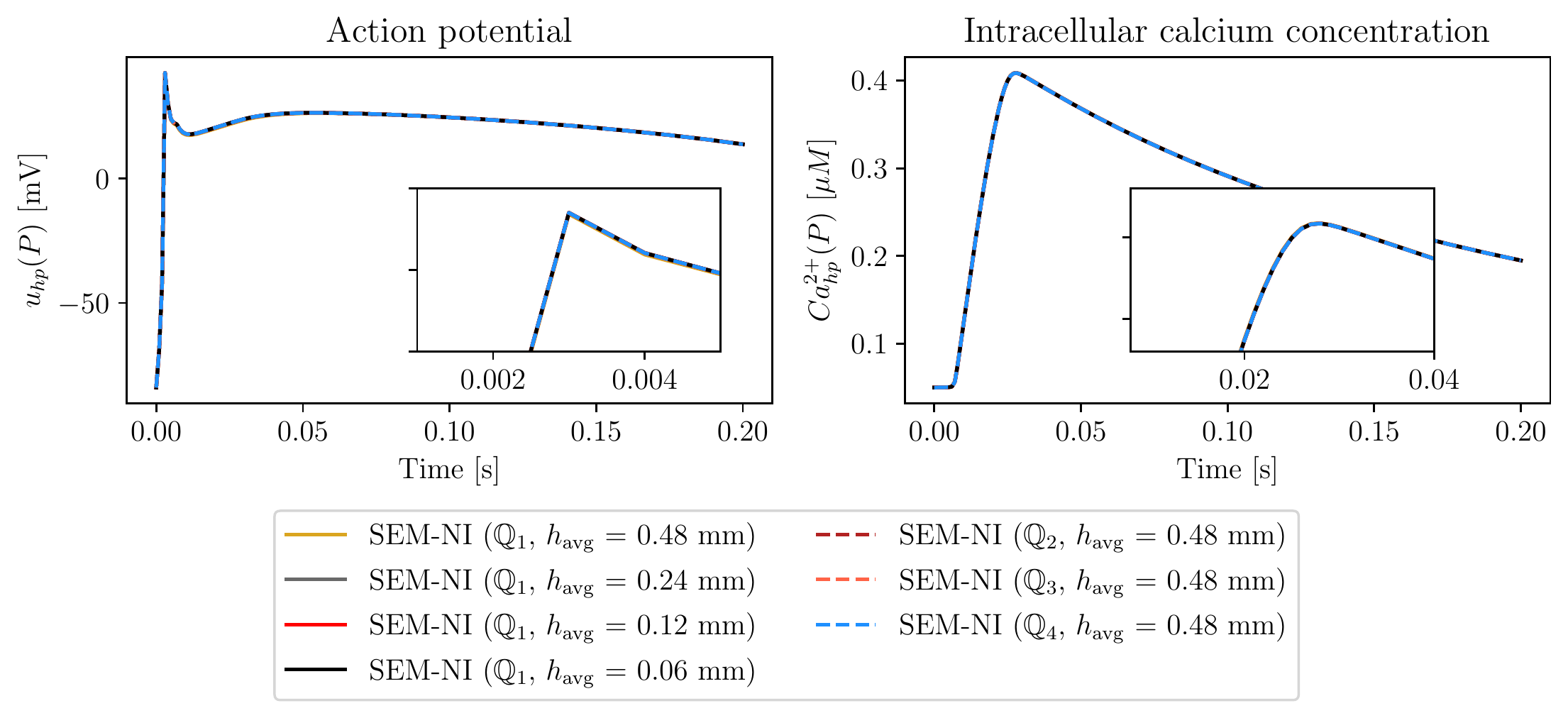}
       \end{subfigure}
      \caption{\emph{Niederer benchmark}. Minimum (top),
    average (second), maximum
    (third)  and point values (bottom) of the action
    potential $u$ and intracellular calcium concentration $Ca^{2+}$ over time for a
    slab of cardiac tissue. $P$ is a random point within the computational domain
    away from the initial stimulus. We consider different $hp$ combinations:
    SEM-NI ${\mathbb Q}_1$ to ${\mathbb Q}_4$ and
    $h_\mathrm{avg} = \SI{0.48}{\milli\meter}$ to $h_\mathrm{avg} =
    \SI{0.06}{\milli\meter}$  ($h_\mathrm{avg}$ is the average mesh size).}
      \label{fig: slab_ap_calcium_sem-ni}
    \end{figure}

    In Figures \ref{fig: slab_ap_calcium_sem} and \ref{fig: slab_ap_calcium_sem-ni}
    we show the action potential and the calcium concentration computed with SEM
    and SEM--NI, respectively, over time. More precisely, the minimum,
    average, maximum and point values are plotted, where the $\max$, $\min$, and $\mathrm{mean}$ functions are evaluated on the set of nodes of the mesh.
    We notice that the convergence
    is faster for increasing $p$ rather than for vanishing
    $h$.

    At each node ${\bf x}_i$ of the mesh, we also compute the activation time \(\tau\) as the time instant when the approximation of the transmembrane potential $u$ exhibits maximum derivative, \textit{i.e.}

    \begin{equation}
      \label{eq:tau}
      \tau\left(\mathbf{x}_i\right) = \argmax_{t} \left|\frac{\partial u}{\partial t}\left(\mathbf{x}_i, t\right)\right|.
    \end{equation}

    \noindent In the formula above $t$ spans over the discrete set of time steps and the time derivative is approximated via the same scheme used for the time discretization of problem \eqref{monodomain}.

\begin{figure}
        \centering
        \begin{subfigure}[t]{0.9\linewidth}
            \centering
            \includegraphics[width=\textwidth]{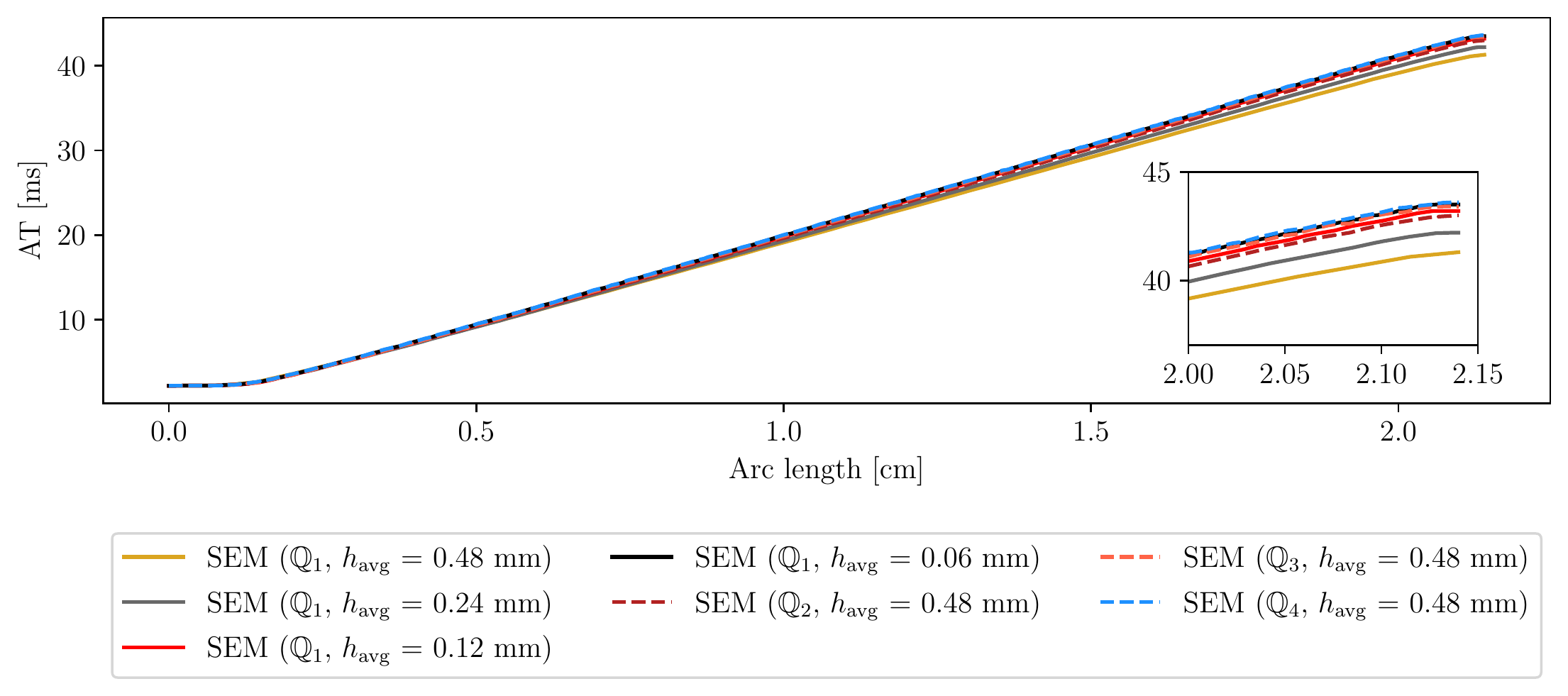}
            \caption{SEM.}
        \end{subfigure}
        \\
        \begin{subfigure}[t]{0.9\linewidth}
            \centering
            \includegraphics[width=\textwidth]{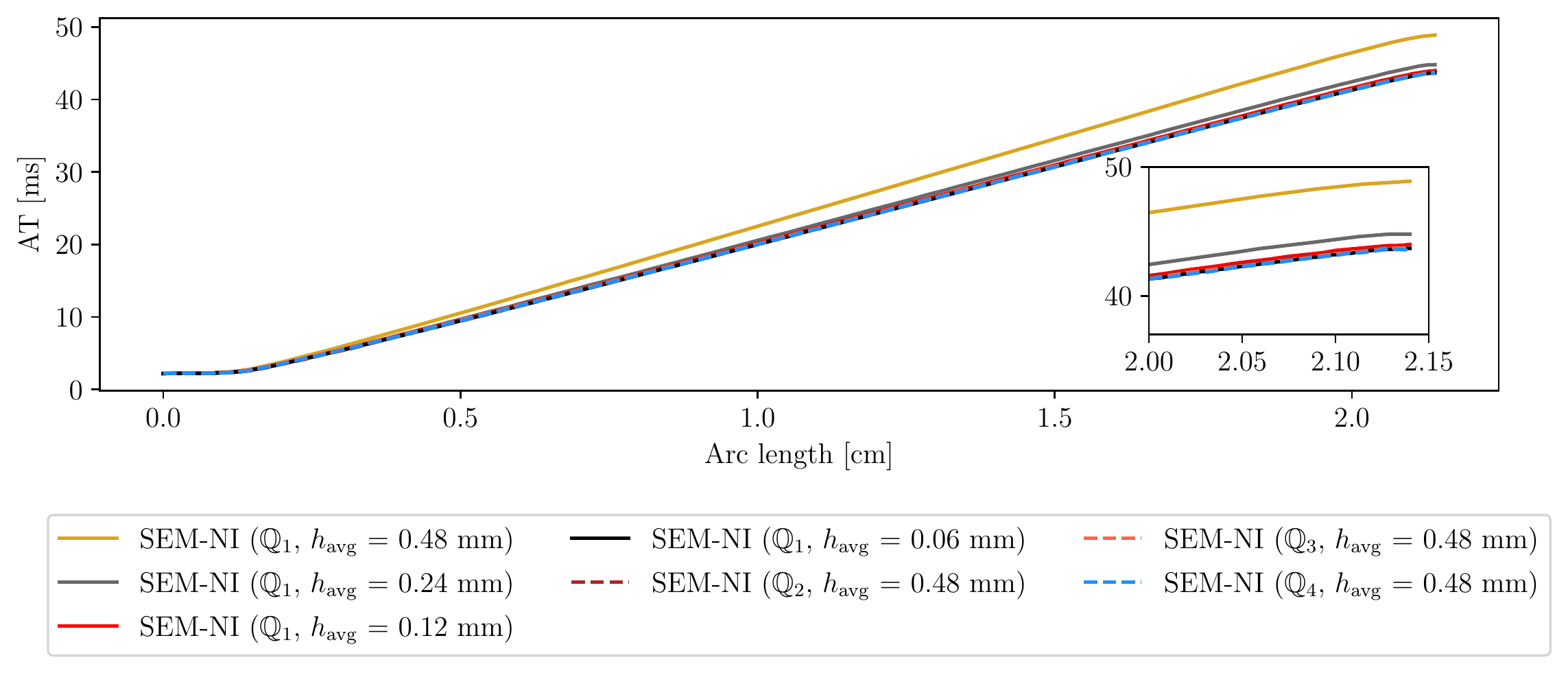}
            \caption{SEM--NI.}
        \end{subfigure}
        \caption{\emph{Niederer benchmark}. Activation times computed along the
    diagonal of the slab (see Figure~\ref{fig:slab}), with different choices of the
    local space (${\mathbb Q}_1$ to ${\mathbb Q}_4$) and mesh refinements
    ($h_\mathrm{avg} = \SI{0.48}{\milli\meter}$ to $h_\mathrm{avg} =
    \SI{0.06}{\milli\meter}$).}
        \label{fig:slab_activation_times}
    \end{figure}

    In Figure~\ref{fig:slab_activation_times} we show the activation times along the
    slab diagonal, for different choices of the local space (from ${\mathbb Q}_1$
    to ${\mathbb Q}_4$) and mesh refinements. As the error accumulates over the
    diagonal, the inset plots show a zoom around the right endpoint. Such results
    demonstrate that high polynomial degrees $p$, even with a coarse mesh size $h$,
    lead to a faster convergence rate compared to the small--$p$, small--$h$ scenario.

\begin{figure}
    \centering
    \begin{subfigure}{.9\linewidth}
    \centering
    \includegraphics[keepaspectratio, width=\textwidth]{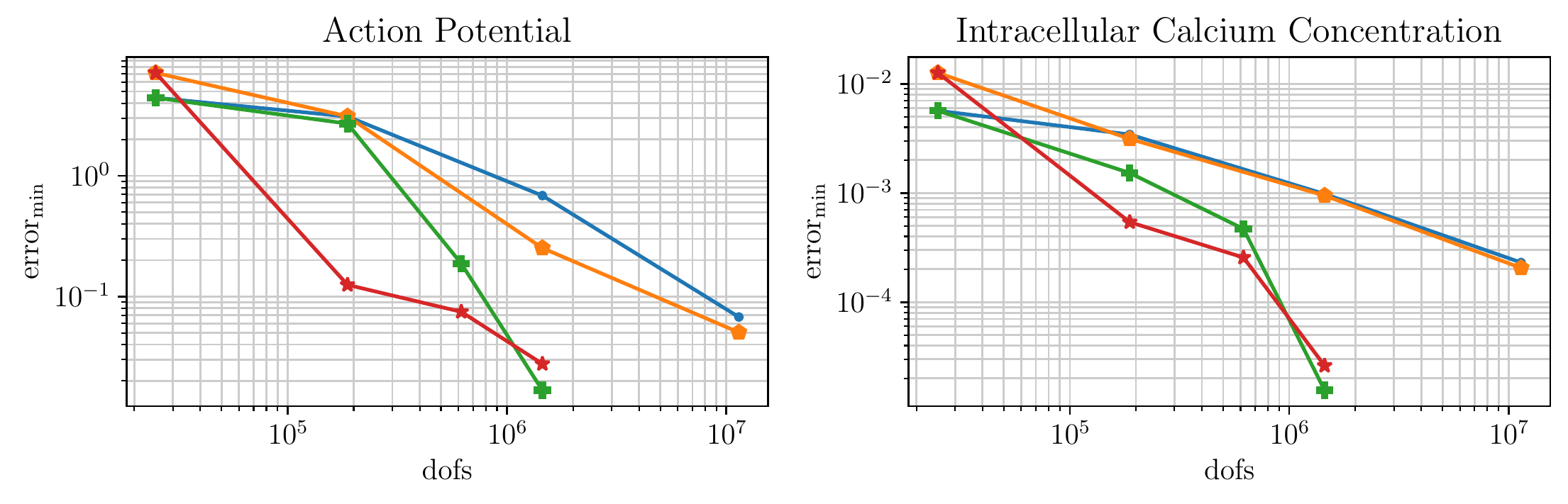}
    \end{subfigure}
    \hfill
    \begin{subfigure}{.9\linewidth}
    \centering
    \includegraphics[keepaspectratio, width=\textwidth]{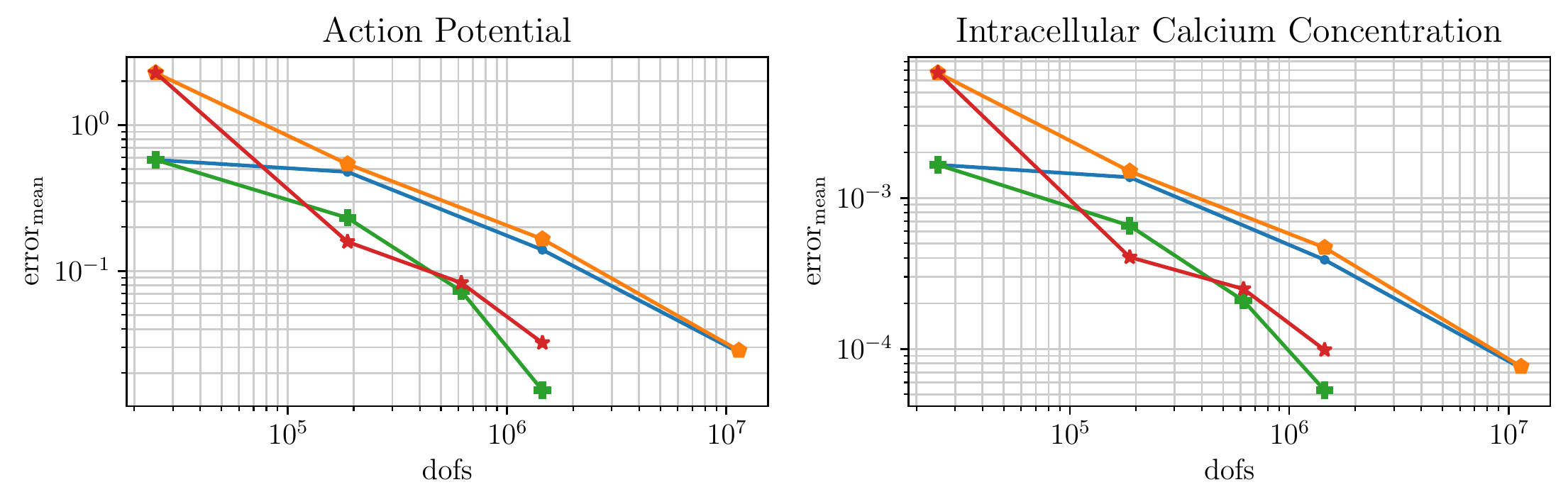}
    \end{subfigure}
    \begin{subfigure}{.9\linewidth}
    \centering
    \includegraphics[keepaspectratio, width=\textwidth]{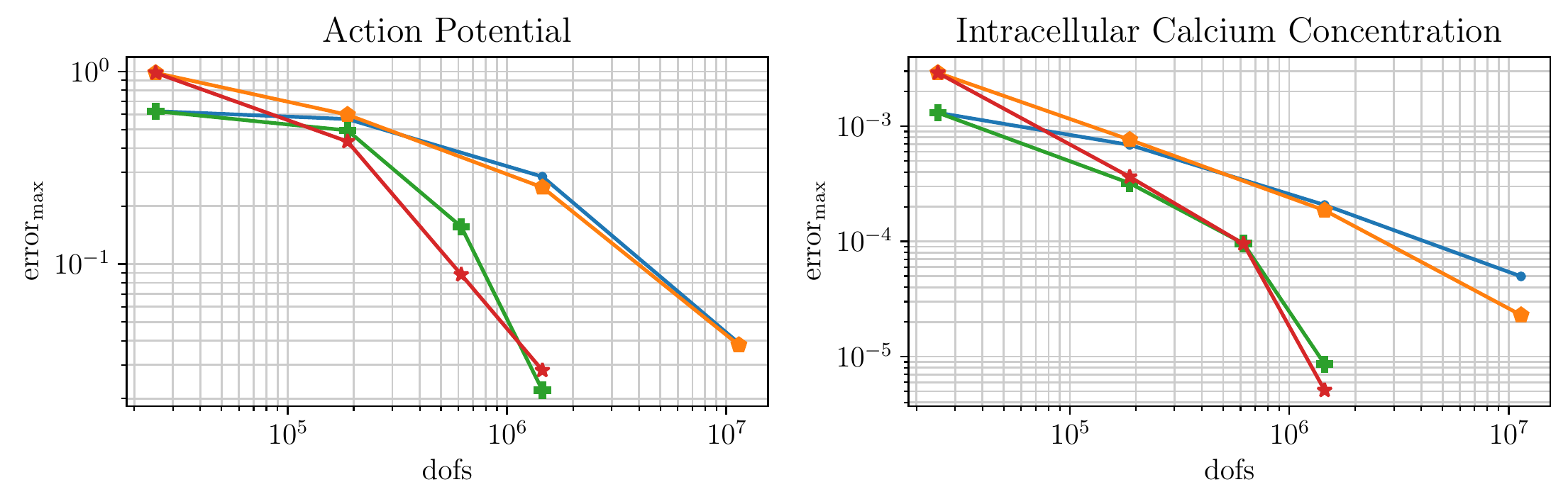}
    \end{subfigure}
    \hfill
    \begin{subfigure}{.9\linewidth}
    \centering
    \includegraphics[keepaspectratio, width=\textwidth]{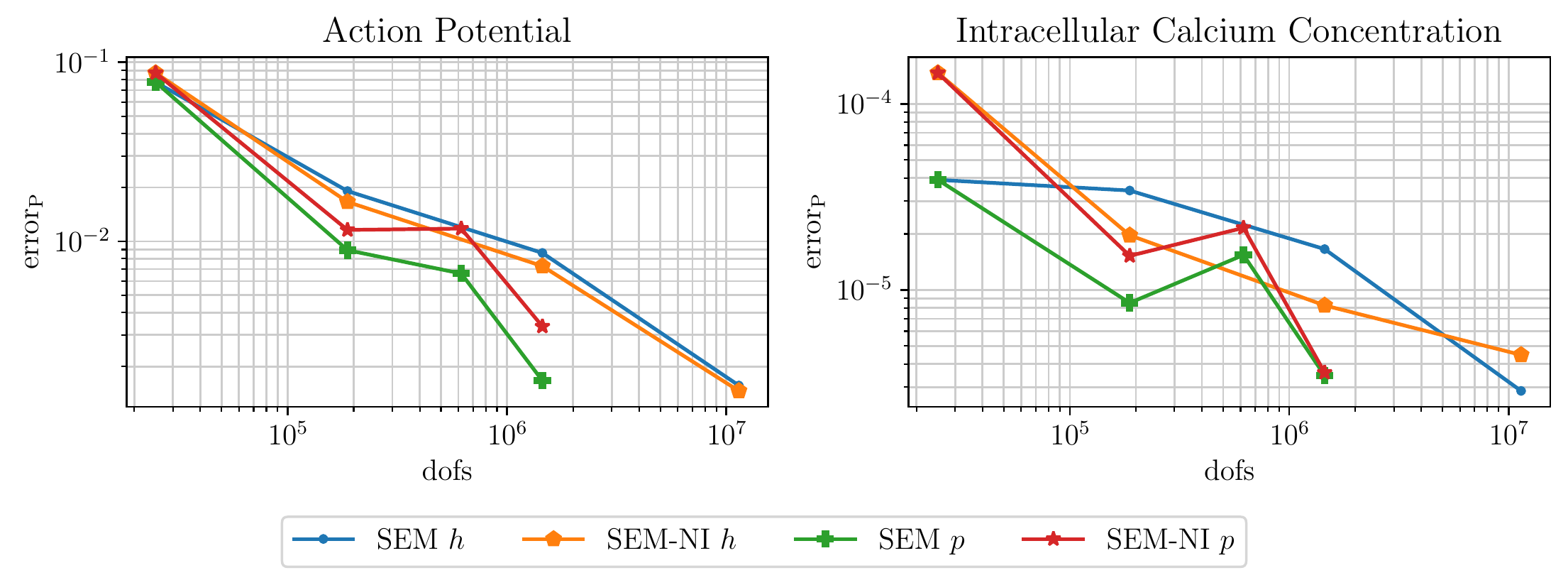}
    \end{subfigure}
    \caption{\emph{Niederer benchmark}. Errors
    for action potential $u$ and intracellular calcium concentration $Ca^{2+}$
    versus the number of mesh points (DOFs) used in
    a slab of cardiac tissue.}
      \label{fig:slab_errors_vs_dofs}
    \end{figure}

    To better investigate the comparison between SEM and SEM--NI,
    in Figure~\ref{fig:slab_errors_vs_dofs} we show the quantities

    \begin{subequations}
    \label{eq:errors}
    \begin{align}
    \mathrm{err}_{\max}&=\left(\Delta t\sum_n \left| \max_{{\bf x}} u_{hp}({\bf
    x},t_n)-\max_{{\bf x}} u_\mathrm{ref}({\bf x},t_n)\right|^2\right)^{1/2} \\
    \mathrm{err}_{\min}&=\left(\Delta t\sum_n \left| \min_{{\bf x}} u_{hp}({\bf
    x},t_n)-\min_{{\bf x}} u_\mathrm{ref}({\bf x},t_n)\right|^2\right)^{1/2} \\
    \mathrm{err}_{\mean}&=\left(\Delta t\sum_n \left| \mean_{{\bf x}}
    u_{hp}({\bf x},t_n)-
    \mean_{{\bf x}}u_\mathrm{ref}({\bf x},t_n)\right|^2\right)^{1/2}\\
    \mathrm{err}_P&=\left(\Delta t\sum_n \left| u_{hp}(P,t_n)- u_\mathrm{ref}(P,t_n)\right|^2\right)^{1/2}
    \end{align}
    \end{subequations}
    versus the total number of mesh points, for both the fully discrete SEM and SEM--NI solutions
    $u_{hp}$.
    Our reference solution $u_\mathrm{ref}$  has been computed with ${\mathbb
    Q}_4-$SEM on a grid with average mesh size $h_\mathrm{avg} = \SI{0.24}{\milli\meter}$,
    for a total of \(11'401'089\) mesh points.
    $P$ is a random point within the computational domain away from the initial stimulus.
    The $\max$, $\min$, and $\mathrm{mean}$ functions are evaluated on the set of nodes of the mesh.
    The number of mesh points increases by reducing $h$ for both ``SEM $h$''  and ``SEM--NI $h$'',
    while it increases with $p$ for both ``SEM $p$'' and ``SEM--NI $p$''. The numerical results confirm the typical behaviour of SEM and SEM--NI discretizations,
    \textit{i.e.} the errors decrease faster by increasing $p$ rather than by decreasing $h$.
    Moreover, we notice that SEM and SEM--NI errors behave quite similarly, with a
    slight advantage for SEM.

\begin{figure}
    \centering
    \begin{subfigure}{.9\linewidth}
    \centering
    \includegraphics[keepaspectratio, width=\textwidth]{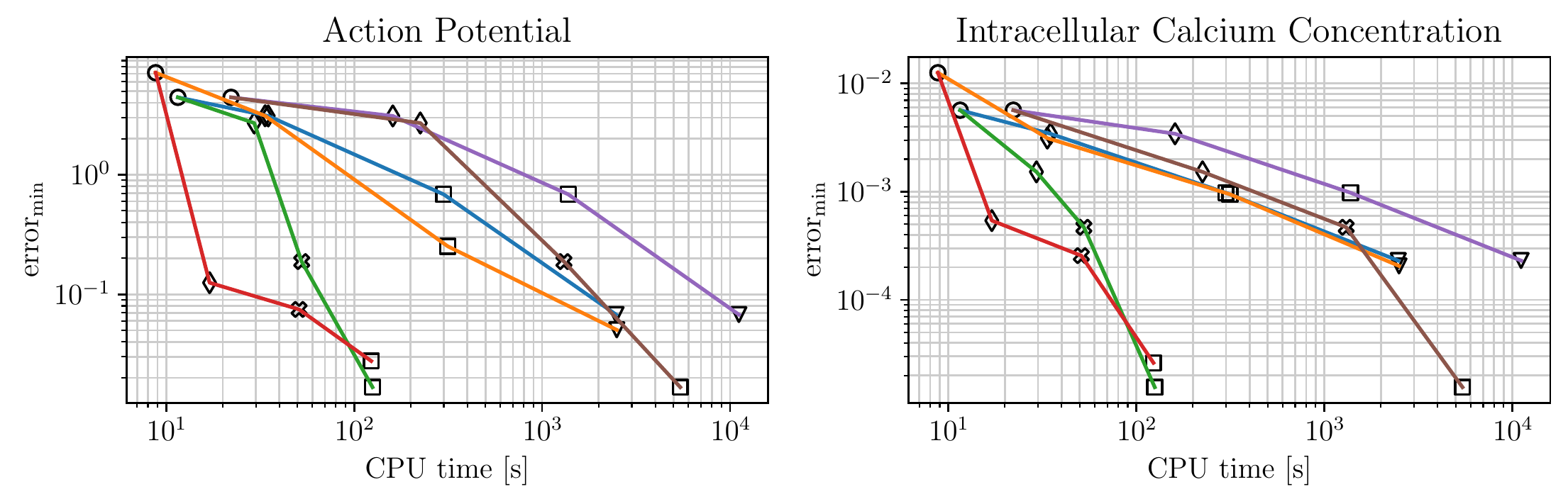}
    \end{subfigure}
    \hfill
    \begin{subfigure}{.9\linewidth}
    \centering
    \includegraphics[keepaspectratio, width=\textwidth]{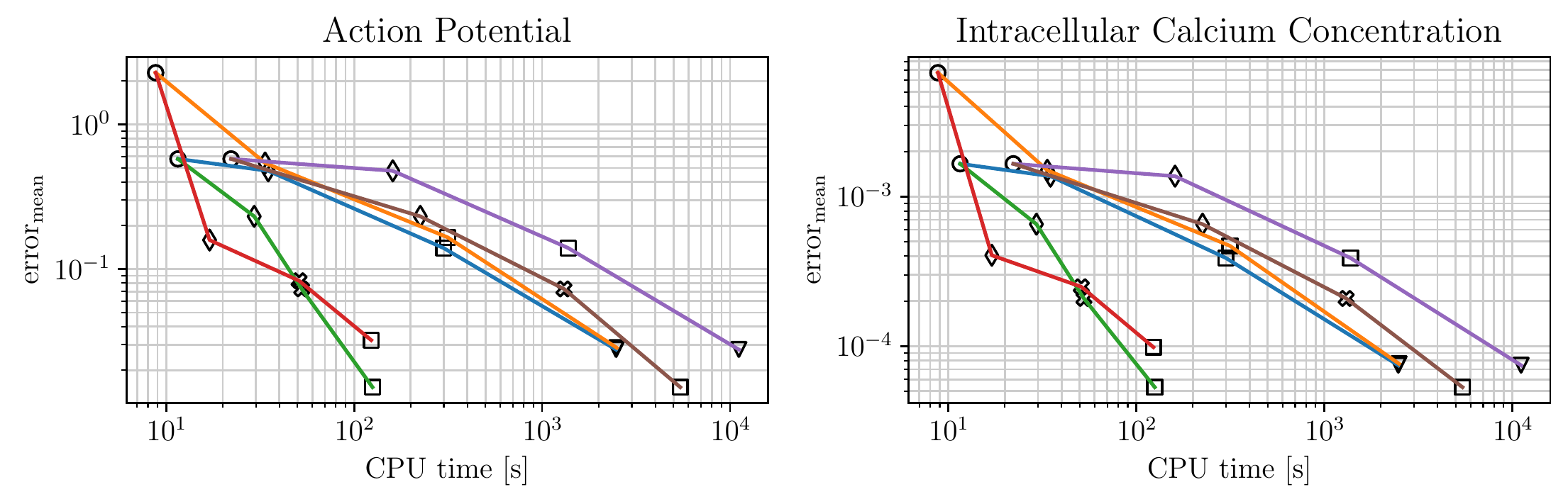}
    \end{subfigure}
    \begin{subfigure}{.9\linewidth}
    \centering
    \includegraphics[keepaspectratio, width=\textwidth]{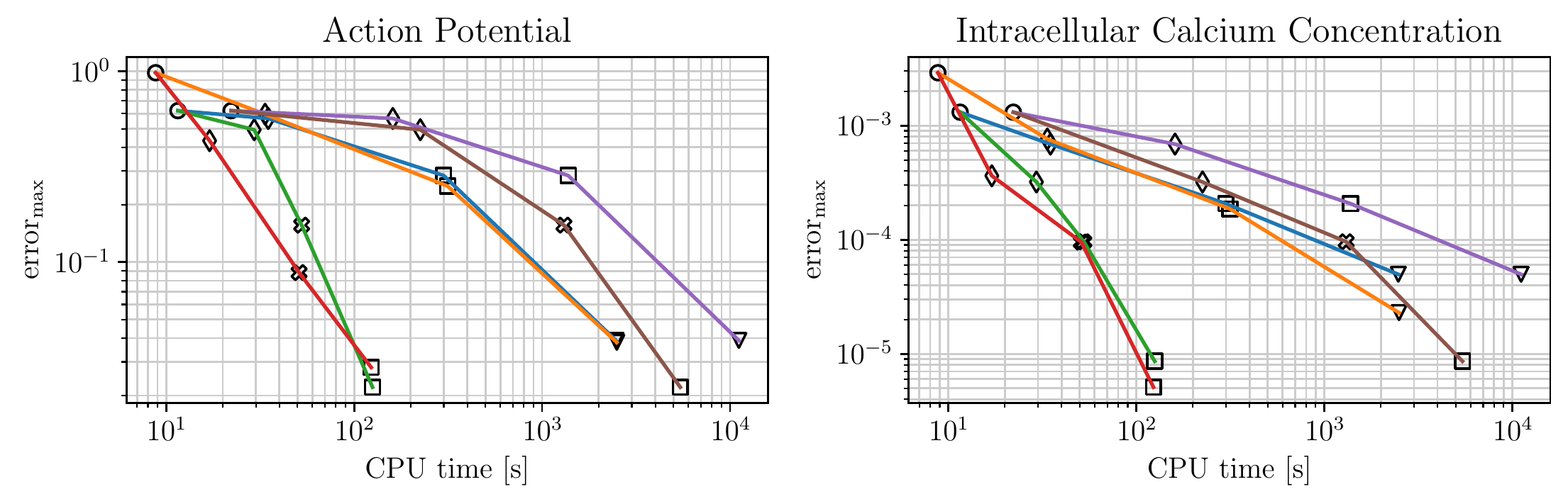}
    \end{subfigure}
    \hfill
    \begin{subfigure}{.9\linewidth}
    \centering
    \includegraphics[keepaspectratio, width=\textwidth]{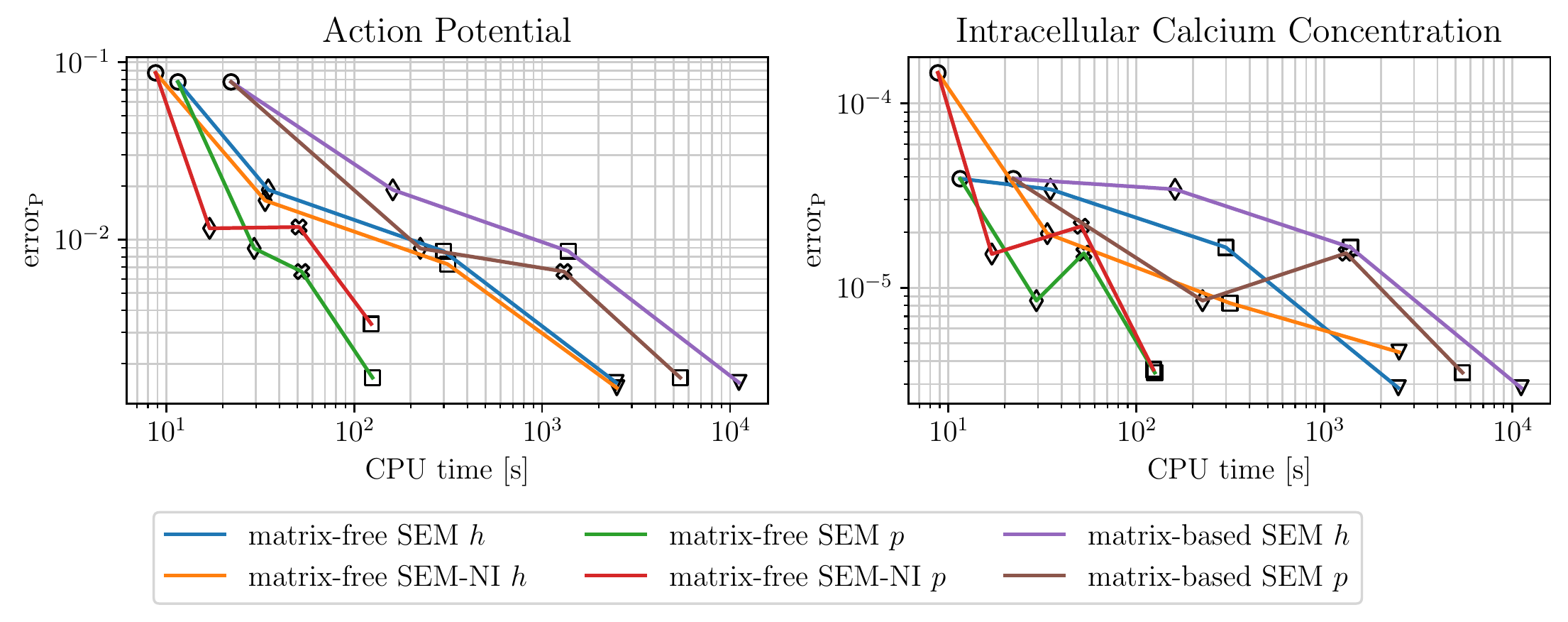}
    \end{subfigure}
    \caption{\emph{Niederer benchmark}. Errors of action potential $u$ and
    intracellular calcium concentration $Ca^{2+}$
    versus CPU time required to solve all the linear systems in
    a slab of cardiac tissue. Equal symbols identify the same number of mesh
    points: $\circ$ $25'025$, $\diamondsuit$ $187'425$, $\times$ $618'529$, $\Box$
    $1'449'665$, $\triangledown$ $11'401'089$.}
      \label{fig:slab_errors_vs_times}
    \end{figure}

    In Tables \ref{tab: SEM_SEM--NI_p}--\ref{tab: matrix--free_matrix--based_h} we
    report the CPU time required by the linear solver for the whole numerical simulation (the times are cumulative over all time
    steps), for SEM and SEM--NI discretizations, matrix--free and matrix--based solvers. We refer to Section~\ref{sec:mf_vs_mb} for the details of the different algorithmic phases.
    Furthermore, in Figure~\ref{fig:slab_errors_vs_times} we plot the errors \eqref{eq:errors} versus the CPU time required to solve all the linear systems along
    the whole numerical simulation. For SEM--NI we only report the times relative to the
    matrix--free solver, while for SEM we report the times for both the
    matrix--free and matrix--based solvers. The same symbol (circle, square,
    diamond, and cross) refers to the numerical simulations carried out on the meshes with the same
    number of DOFs.
    If we compare the errors and the CPU--times of ${\mathbb Q}_1-$SEM, $h_{\mathrm{avg}} = \SI{0.12}{\milli\meter}$
    with those of
    ${\mathbb Q}_4-$SEM, $h_{\mathrm{avg}} = \SI{0.48}{\milli\meter}$ (these two
    configurations share the same number \(1'449'665\) of mesh nodes),
    we notice that the errors of ${\mathbb Q}_4-$SEM are at most about 1/3 -- 1/2 of that
    of ${\mathbb Q}_1-$SEM and the ratio between the corresponding CPU times is about
    40\%. Thus, we conclude that ${\mathbb Q}_4-$SEM outperforms ${\mathbb
    Q}_1-$SEM (that is ${\mathbb Q}_1-$FEM).

    \begin{table}
        \hspace*{-.7cm}
        {\fontsize{8}{9}\selectfont{
        \begin{tabular}{ | c c c | c c c c |}
            \toprule
            \makecell{Mesh points \\ number } & \makecell{Cells \\ number} &
    \makecell{Local \\ space} & \makecell{Linear solver \\ SEM [\si{\second}]} & \makecell{Linear solver \\ SEM--NI [\si{\second}]} & \makecell{Assemble rhs \\ SEM [\si{\second}]} & \makecell{Assemble rhs \\ SEM--NI [\si{\second}]} \\
            \midrule
            $25'025$    & $21'888$ & $\mathbb{Q}_1$ & 10.769 & 8.031  & 0.784  & 0.766  \\
            $187'425$   & $21'888$ & $\mathbb{Q}_2$ & 14.694 & 12.383 & 4.710  & 4.645  \\
            $618'529$   & $21'888$ & $\mathbb{Q}_3$ & 37.733 & 36.419 & 14.867 & 14.542 \\
            $1'449'665$ & $21'888$ & $\mathbb{Q}_4$ & 91.380 & 90.370 & 33.899 & 32.920 \\
            \bottomrule
        \end{tabular}
    }}
        \caption{\emph{Niederer benchmark}. Computational times for SEM and SEM--NI, with a fixed average mesh size $h_\mathrm{avg} = \SI{0.48}{\milli\meter}$ and $p$ ranging from 1 to 4. Matrix--free solver.}
        \label{tab: SEM_SEM--NI_p}
    \end{table}

    \begin{table}
        \hspace*{-.8cm}
        {\fontsize{8}{9}\selectfont{
        \begin{tabular}{ | c c c | c c c c |}
            \toprule
            \makecell{Mesh points \\ number }  & \makecell{Cells \\ number} & \makecell{$h_\mathrm{avg}$ \\ \text{[\si{\milli\meter}]}} & \makecell{Linear solver \\ SEM [\si{\second}]} & \makecell{Linear solver \\ SEM--NI [\si{\second}]} & \makecell{Assemble rhs \\ SEM [\si{\second}]} & \makecell{Assemble rhs \\ SEM--NI [\si{\second}]} \\
            \midrule
            $25'025$    & $21'888$     & 0.48 & 10.769   & 8.031    & 0.784   & 0.766   \\
            $187'425$   & $175'104$    & 0.24 & 29.157   & 27.951   & 5.771   & 5.656   \\
            $1'449'665$ & $1'400'832$  & 0.12 & 256.295  & 270.959  & 42.783  & 43.548  \\
            $11'401'089$ & $11'206'656$ & 0.06 & 2137.329 & 2158.751 & 336.272 & 336.641 \\
            \bottomrule
        \end{tabular}
    }}
        \caption{\emph{Niederer benchmark}. Computational times for SEM and SEM--NI, with an average mesh size $h_\mathrm{avg}$ ranging from $\SI{0.48}{\milli\meter}$ to $\SI{0.06}{\milli\meter}$ and $p$ = 1. Matrix--free solver.}
        \label{tab: SEM_SEM--NI_h}
    \end{table}

    \begin{table}
            \hspace*{-1.2cm}
        {\fontsize{8}{9}\selectfont{
        \begin{tabular}{ | c c c | c c c c |}
            \toprule
            \makecell{Mesh points \\ number }  & \makecell{Cells \\ number} &
    \makecell{Local \\ space} & \makecell{Linear solver \\ matrix--free [\si{\second}]} & \makecell{Assembly phase \\ matrix--free [\si{\second}]} & \makecell{Linear solver \\ matrix--based [\si{\second}]} & \makecell{Assembly phase \\ matrix--based [\si{\second}]} \\
            \midrule
            $25'025$    & $21'888$ & $\mathbb{Q}_1$ & 10.769 & 0.784  & 4.086    & 18.076   \\
            $187'425$   & $21'888$ & $\mathbb{Q}_2$ & 14.694 & 4.710  & 44.200   & 180.705  \\
            $618'529$   & $21'888$ & $\mathbb{Q}_3$ & 37.733 & 14.867 & 343.243  & 963.549  \\
            $1'449'665$ & $21'888$ & $\mathbb{Q}_4$ & 91.380 & 33.899 & 1557.144 & 3874.602 \\
            \bottomrule
        \end{tabular}
    }}
        \caption{\emph{Niederer benchmark}. Computational times for matrix--free and matrix--based solvers, SEM ${\mathbb Q}_p$ with a fixed average mesh size $h_\mathrm{avg} = \SI{0.48}{\milli\meter}$ and $p$ ranging from 1 to 4.}
        \label{tab: matrix--free_matrix--based_p}
    \end{table}

    \begin{table}
        \hspace*{-1.4cm}
        {\fontsize{8}{9}\selectfont{
        \begin{tabular}{ | c c c | c c c c |}
            \toprule
            \makecell{Mesh points \\ number }  & \makecell{Cells \\ number} & \makecell{$h_\mathrm{avg}$ \\ \text{[\si{\milli\meter}]}} & \makecell{Linear solver \\ matrix--free [\si{\second}]} & \makecell{Assembly phase \\ matrix--free [\si{\second}]} & \makecell{Linear solver \\ matrix--based [\si{\second}]} & \makecell{Assembly phase \\ matrix--based [\si{\second}]} \\
            \midrule
            $25'025$    & $21'888$     & 0.48 & 10.769   & 0.784   & 4.086    & 18.076   \\
            $187'425$   & $175'104$    & 0.24 & 29.157   & 5.771   & 15.373   & 145.244  \\
            $1'449'665$ & $1'400'832$  & 0.12 & 256.295  & 42.783  & 204.809  & 1171.724 \\
            $11'401'089$ & $11'206'656$ & 0.06 & 2137.329 & 336.272 & 1867.746 & 9266.635 \\
            \bottomrule
        \end{tabular}
    }}
        \caption{\emph{Niederer benchmark}. Computational times for matrix--free and
    matrix--based, SEM ${\mathbb Q}_1$ with an average mesh size $h_\mathrm{avg}$ ranging from $\SI{0.48}{\milli\meter}$ to $\SI{0.06}{\milli\meter}$.}
        \label{tab: matrix--free_matrix--based_h}
    \end{table}

    For the comparison between matrix--free and matrix--based
    solvers, we notice that the former one is always faster, and the gain
    of matrix--free over matrix--based solver increases with the polynomial degree $p$.
    More precisely, the speed--up factors
    are shown in Table \ref{tab:speed-up-p} when $h_{\mathrm{avg}} =
    \SI{0.48}{\milli\meter}$, and in Table \ref{tab:speed-up-h} when $p=1$.

    \begin{table}
    \begin{center}
    \begin{tabular}{| r  | c c c c |}
    \toprule
    Local space &  ${\mathbb Q}_1$ &  ${\mathbb
    Q}_2$ & ${\mathbb Q}_3$ & ${\mathbb Q}_4$\\
    \midrule
    $\displaystyle \frac{\mathrm{(CPU\ time)}_{\mathrm {mb}}} {\mathrm{(CPU\ time)}_{\mathrm{mf}}}$
    & $\sim 2$ & $\sim 12$ & $\sim 25$ & $\sim 45$\\
    \bottomrule
    \end{tabular}
    \end{center}
    \caption{\emph{Niederer benchmark}. Speed--up of the matrix--free solver over
    the matrix--based one when $h_{\mathrm{avg}} =
    \SI{0.48}{\milli\meter}$.}\label{tab:speed-up-p}
    \end{table}

    \begin{table}
    \begin{center}
    \begin{tabular}{| r  | c c c c |}
    \toprule
    $h_{\mathrm{avg}}$ &  0.48 & 0.24 & 0.12 & 0.06\\
    \midrule
    $\displaystyle \frac{\mathrm{(CPU\ time)}_{\mathrm {mb}}} {\mathrm{(CPU\ time)}_{\mathrm{mf}}}$
    & $\sim 2$ & $\sim 5$ & $\sim 5$ & $\sim 5$\\
    \bottomrule
    \end{tabular}
    \end{center}
    \caption{\emph{Niederer benchmark}. Speed--up of the matrix--free solver over
    the matrix--based one when $p=1$.}
    \label{tab:speed-up-h}
    \end{table}

    Moreover, from Table \ref{tab: matrix--free_matrix--based_percentages} we
    observe that, in a matrix--based electrophysiological simulation,
    most of the computational time is spent to solve the linear system associated
    with the monodomain equation.
    On the contrary, in the matrix--free solver most of the computational time is
    devoted to the ionic model. This means that the cost for solving the linear
    system has been highly optimized.

    \begin{table}
        \hspace*{-1.4cm}
        {\fontsize{8}{9}\selectfont{
        \begin{tabular}{| l | c c c |}
            \toprule
            \makecell{Solver} & \makecell{Monodomain solver} & \makecell{Monodomain assembly} & \makecell{Ionic model solver} \\
            \midrule
            Matrix--based ($\mathbb{Q}_1$, $h_\mathrm{avg} = \SI{0.12}{\milli\meter}$) & 10.54 \% & 60.41 \% & 29.05 \% \\
            Matrix--free  ($\mathbb{Q}_4$, $h_\mathrm{avg} = \SI{0.48}{\milli\meter}$) & 14.95 \% & 3.36  \% & 81.69 \% \\
            \bottomrule
        \end{tabular}
    }}
        \caption{\emph{Niederer benchmark}. Matrix--free and matrix--based percentages for the assembling and solving phases.
        Note that in the matrix--free case we only need to assemble the right--hand
    side vector, as there is no matrix. Moreover, these percentages are computed
    without taking into account all other phases of the numerical simulation, such as mesh allocation, fiber generation and output of the results.}
        \label{tab: matrix--free_matrix--based_percentages}
    \end{table}

    Finally, we compare the performance of the AMG and GMG preconditioners,
    used by the matrix--based and matrix--free solvers, respectively. In Tables
    \ref{tab: matrix--free_matrix--based_CG_p} and \ref{tab:
    matrix--free_matrix--based_CG_h} we show the average number of iterations
    required by the PCG method to solve the linear system \eqref{linsys} for different combinations of \(h\) and \(p\).
    We notice that, for the values of $h$ and $p$ considered here, both the AMG and GMG preconditioners appear to be optimal in the number of PCG iterations versus both $h$ and $p$.
    As a matter of fact, the average number of iterations is about 1.0 (matrix--based) and 1.8 (matrix--free) for all $hp$ configurations. More precisely, the number of iterations throughout all the simulations ranges from 1 to 4.

    \begin{table}
            \hspace*{-1.2cm}
        {\fontsize{8}{9}\selectfont{
        \begin{tabular}{ | c c c | c c c |}
            \toprule
            \makecell{Mesh points \\ number}  & \makecell{Cells \\ number} &
    \makecell{Local \\ space} & \makecell{Matrix--free (SEM)\\ GMG preconditioner} &
    \makecell{Matrix--free (SEM--NI)\\ GMG preconditioner} & \makecell{Matrix--based
    (SEM)\\ AMG preconditioner} \\
            \midrule
            $25'025$    & $21'888$ & $\mathbb{Q}_1$ & 1.6362 & 1.8126 & 0.9780  \\
            $187'425$   & $21'888$ & $\mathbb{Q}_2$ & 1.8056 & 1.8581 & 1.0025  \\
            $618'529$   & $21'888$ & $\mathbb{Q}_3$ & 1.7906 & 1.8161 & 1.0205  \\
            $1'449'665$ & $21'888$ & $\mathbb{Q}_4$ & 1.7371 & 1.7826 & 1.0250  \\
            \bottomrule
        \end{tabular}
    }}
        \caption{\emph{Niederer benchmark}. Average number of CG iterations for
    matrix--free (SEM, SEM--NI) and matrix--based (SEM) solvers, ${\mathbb Q}_p$ with a fixed average mesh size $h_\mathrm{avg} = \SI{0.48}{\milli\meter}$ and $p$ ranging from 1 to 4.}
        \label{tab: matrix--free_matrix--based_CG_p}
    \end{table}

    \begin{table}
        \hspace*{-1.4cm}
        {\fontsize{8}{9}\selectfont{
        \begin{tabular}{ | c c c | c c c |}
            \toprule
            \makecell{Mesh points \\ number }  & \makecell{Cells \\ number} &
    \makecell{$h_\mathrm{avg}$ \\ \text{[\si{\milli\meter}]}} & \makecell{Matrix--free (SEM)\\ GMG
    preconditioner} &
    \makecell{Matrix--free (SEM--NI)\\ GMG preconditioner} & \makecell{Matrix--based (SEM)\\ AMG
    preconditioner} \\
            \midrule
            $25'025$    & $21'888$     & 0.48 & 1.6362 & 1.8126 & 0.9780 \\
            $187'425$   & $175'104$    & 0.24 & 1.5757 & 1.6847 & 0.9855 \\
            $1'449'665$ & $1'400'832$  & 0.12 & 1.4448 & 1.5937 & 1.0060 \\
            $11'401'089$ & $11'206'656$ & 0.06 & 1.4468 & 1.4923 & 1.0180 \\
            \bottomrule
        \end{tabular}
    }}
        \caption{\emph{Niederer benchmark}. Average number of CG iterations for
    matrix--free (SEM, SEM--NI) and matrix--based (SEM) solvers, ${\mathbb Q}_1$ with an average mesh size $h_\mathrm{avg}$ ranging from $\SI{0.48}{\milli\meter}$ to $\SI{0.06}{\milli\meter}$.}
        \label{tab: matrix--free_matrix--based_CG_h}
    \end{table}

    The numerical results shown in this section highlight how much advantageous the matrix--free
    solver with SEM or SEM--NI is for cardiac electrophysiology simulations, with
    respect to the matrix--based solver with low--order FEM.

    Since the matrix--free implementation outperforms the matrix--based one, while
    SEM and SEM--NI provide comparable results in terms of accuracy and efficiency,
    we will employ the matrix--free solver with just the SEM formulation for the
    numerical simulations that we are going to present in the next sections.

    \subsection{Left ventricle}
    \label{sec:lv}

    We report the results for the electrophysiological simulations performed with the Zygote left
    ventricle geometry \cite{Zygote}. The settings of this test case are summarized at the beginning of Section~ \ref{sec:numres}.
    We consider a mesh with $h_{\mathrm{avg}} = \SI{2.0}{\milli\meter}$ and polynomial degree $p$ from 1 to 4. In all cases, we keep the mesh boundary fixed to the one resulting from a linear mapping to neglect the impact of boundary deformation on the accuracy of the numerical simulations.

    In Table \ref{tab: matrix-free_CG_p} we report the number of mesh nodes, the
    number of cells and the average number of iterations required by the PCG method
    to solve the linear system (\ref{linsys}).
    As for the Niederer benchmark, the GMG preconditioner
    turns out to be optimal also for these numerical simulations.
    Indeed, the number of PCG iterations is about 2 along the whole time history
    for any polynomial degree between 1 and 4.

    \begin{table}
            \centering
        {\fontsize{8}{9}\selectfont{
        \begin{tabular}{ | c c c | c |}
            \toprule
            \makecell{Mesh points \\ number}  & \makecell{Cells \\ number} &
    \makecell{Local \\ space} & \makecell{PCG iterations\\ GMG preconditioner} \\
            \midrule
            $159'149$   & $139'684$ & $\mathbb{Q}_1$ & 2.0770 \\
            $1'172'919$ & $139'684$ & $\mathbb{Q}_2$ & 1.9628 \\
            $3'879'415$ & $139'684$ & $\mathbb{Q}_3$ & 1.9455 \\
            $9'116'741$ & $139'684$ & $\mathbb{Q}_4$ & 1.9440 \\
            \bottomrule
        \end{tabular}
    }}
        \caption{\emph{Zygote left ventricle}. Average number of PCG iterations for
    the matrix--free solver with SEM discretization, ${\mathbb Q}_p$ with a fixed average mesh size
    $h_\mathrm{avg} = \SI{2.0}{\milli\meter}$ and $p$ ranging from 1 to 4.}
        \label{tab: matrix-free_CG_p}
    \end{table}

\begin{figure}
        \centering
        \captionsetup{justification=centering,margin=2cm}
        \begin{subfigure}[t]{\linewidth}
            \centering
            \includegraphics[width=\textwidth]{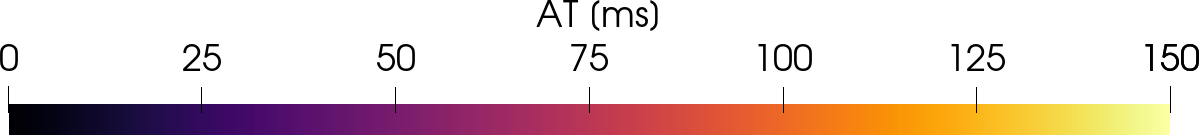}
        \end{subfigure}
        \\[\baselineskip]
        \begin{subfigure}[t]{0.2\linewidth}
            \centering
            \includegraphics[width=\textwidth]{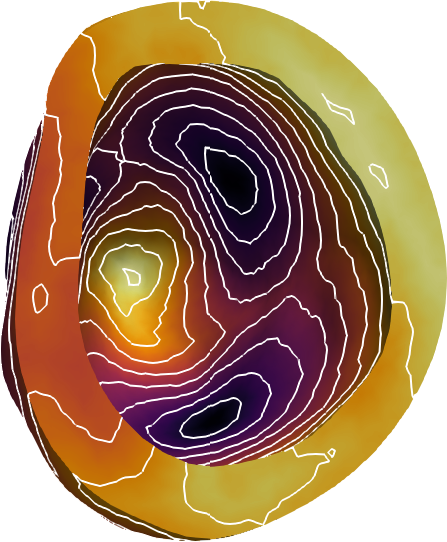}
            \caption*{}
        \end{subfigure}
        \hfill
        \begin{subfigure}[t]{0.2\linewidth}
            \centering
            \includegraphics[width=\textwidth]{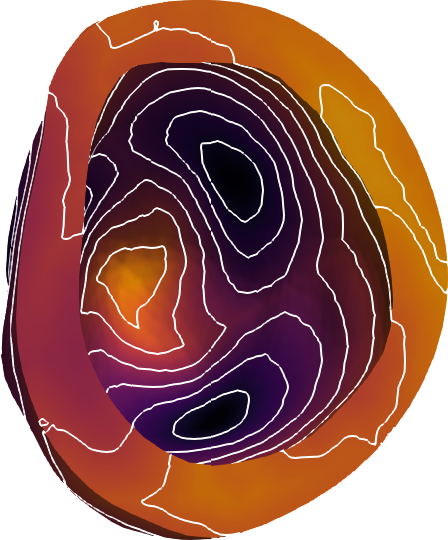}
            \caption*{}
        \end{subfigure}
        \hfill
        \begin{subfigure}[t]{0.2\linewidth}
            \centering
            \includegraphics[width=\textwidth]{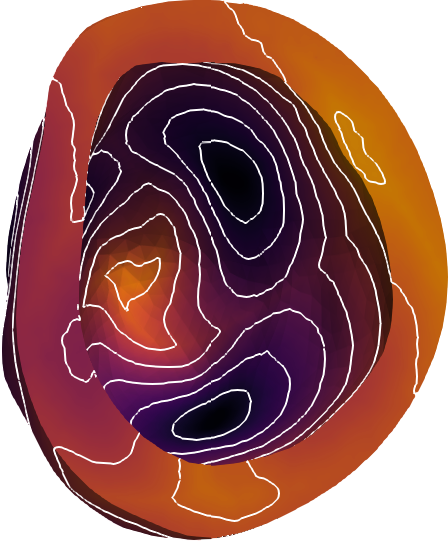}
            \caption*{}
        \end{subfigure}
        \hfill
        \begin{subfigure}[t]{0.2\linewidth}
            \centering
            \includegraphics[width=\textwidth]{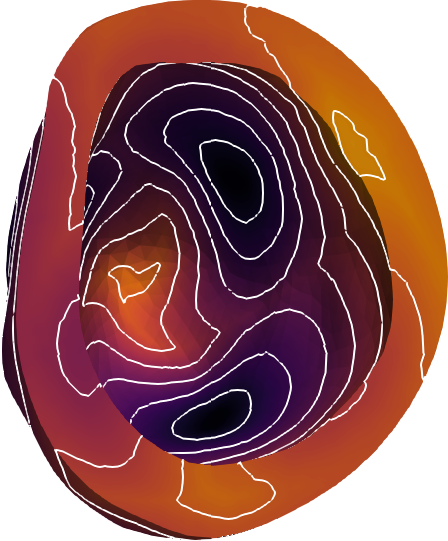}
            \caption*{}
        \end{subfigure}
        \\
        \begin{subfigure}[t]{0.2\linewidth}
            \centering
            \includegraphics[width=\textwidth]{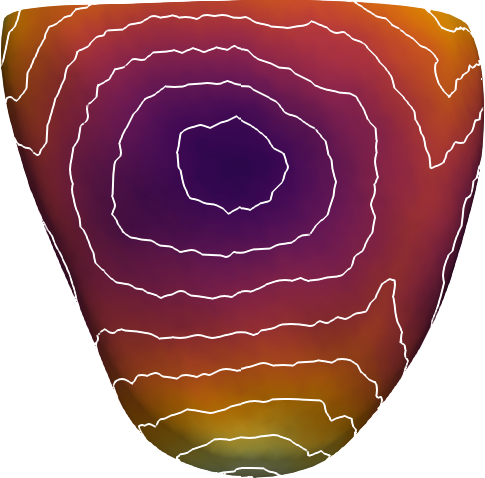}
            \caption*{}
        \end{subfigure}
        \hfill
        \begin{subfigure}[t]{0.2\linewidth}
            \centering
            \includegraphics[width=\textwidth]{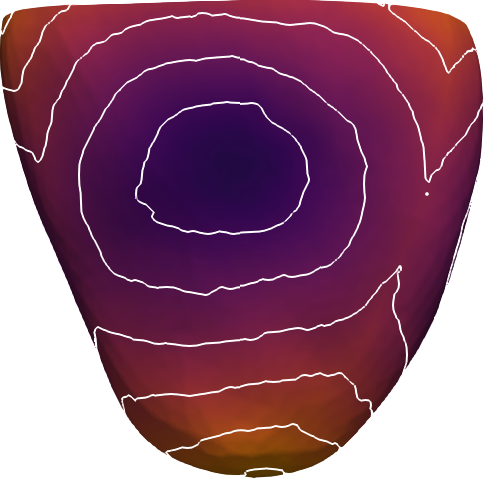}
            \caption*{}
        \end{subfigure}
        \hfill
        \begin{subfigure}[t]{0.2\linewidth}
            \centering
            \includegraphics[width=\textwidth]{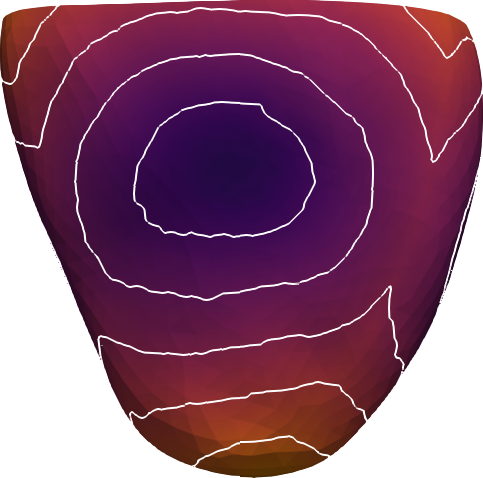}
            \caption*{}
        \end{subfigure}
        \hfill
        \begin{subfigure}[t]{0.2\linewidth}
            \centering
            \includegraphics[width=\textwidth]{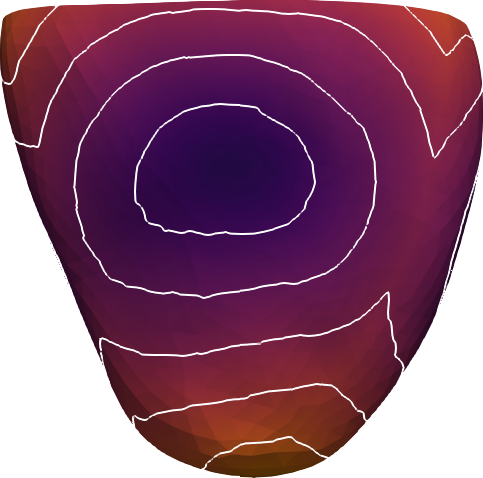}
            \caption*{}
        \end{subfigure}
        \\
        \begin{subfigure}[t]{0.2\linewidth}
            \centering
            \includegraphics[width=\textwidth]{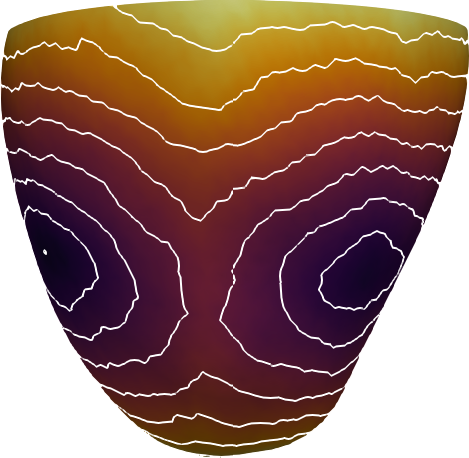}
            \caption*{$159'149 \,\, \text{DOFs}$\\($\mathbb{Q}_1$)}
        \end{subfigure}
        \hfill
        \begin{subfigure}[t]{0.2\linewidth}
            \centering
            \includegraphics[width=\textwidth]{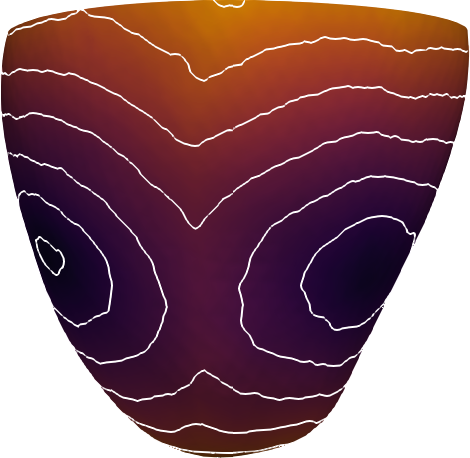}
            \caption*{$1'172'919 \,\, \text{DOFs}$\\($\mathbb{Q}_2$)}
        \end{subfigure}
        \hfill
        \begin{subfigure}[t]{0.2\linewidth}
            \centering
            \includegraphics[width=\textwidth]{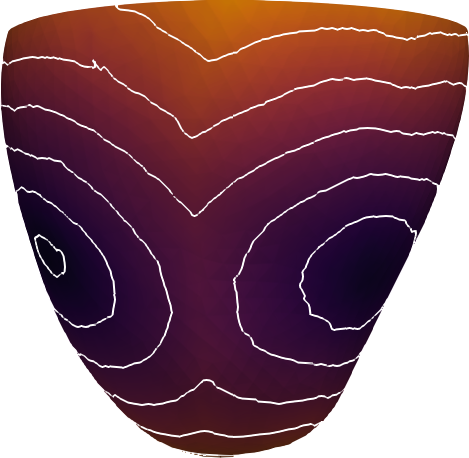}
            \caption*{$3'879'415 \,\, \text{DOFs}$\\($\mathbb{Q}_3$)}
        \end{subfigure}
        \hfill
        \begin{subfigure}[t]{0.2\linewidth}
            \centering
            \includegraphics[width=\textwidth]{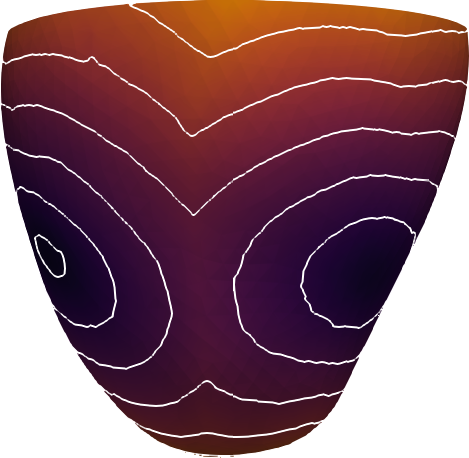}
            \caption*{$9'116'741 \,\, \text{DOFs}$\\($\mathbb{Q}_4$)}
        \end{subfigure}
        \captionsetup{justification=raggedright,margin=0cm}
        \caption{\emph{Zygote left ventricle}. Three views of the activation maps computed with ${\mathbb Q}_p$ elements ($p=1, \dots, 4$) and a fixed average mesh size $h_\mathrm{avg} = \SI{2.0}{\milli\meter}$.}
        \label{fig:activation_time_LV}
    \end{figure}

\begin{figure}
    \centering
    \begin{subfigure}{.95\linewidth}
        \centering
        \includegraphics[keepaspectratio, width=\textwidth]{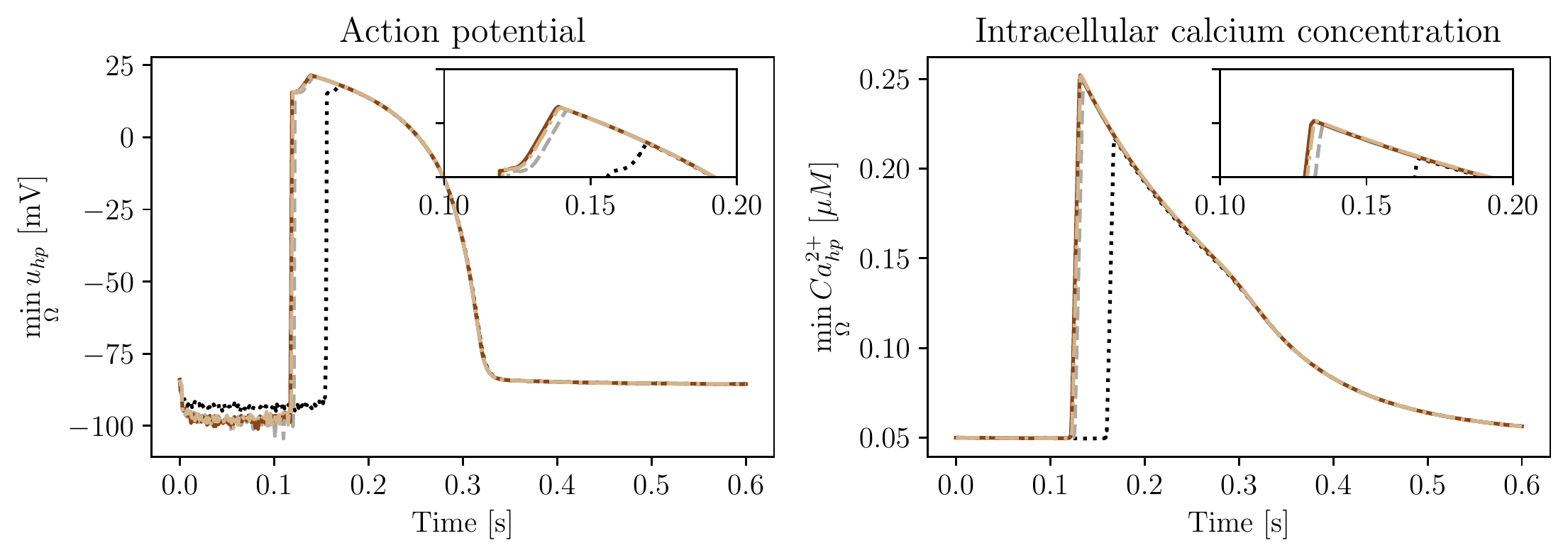}
    \end{subfigure} %

    \hfill

    \begin{subfigure}{.95\linewidth}
        \centering
        \includegraphics[keepaspectratio, width=\textwidth]{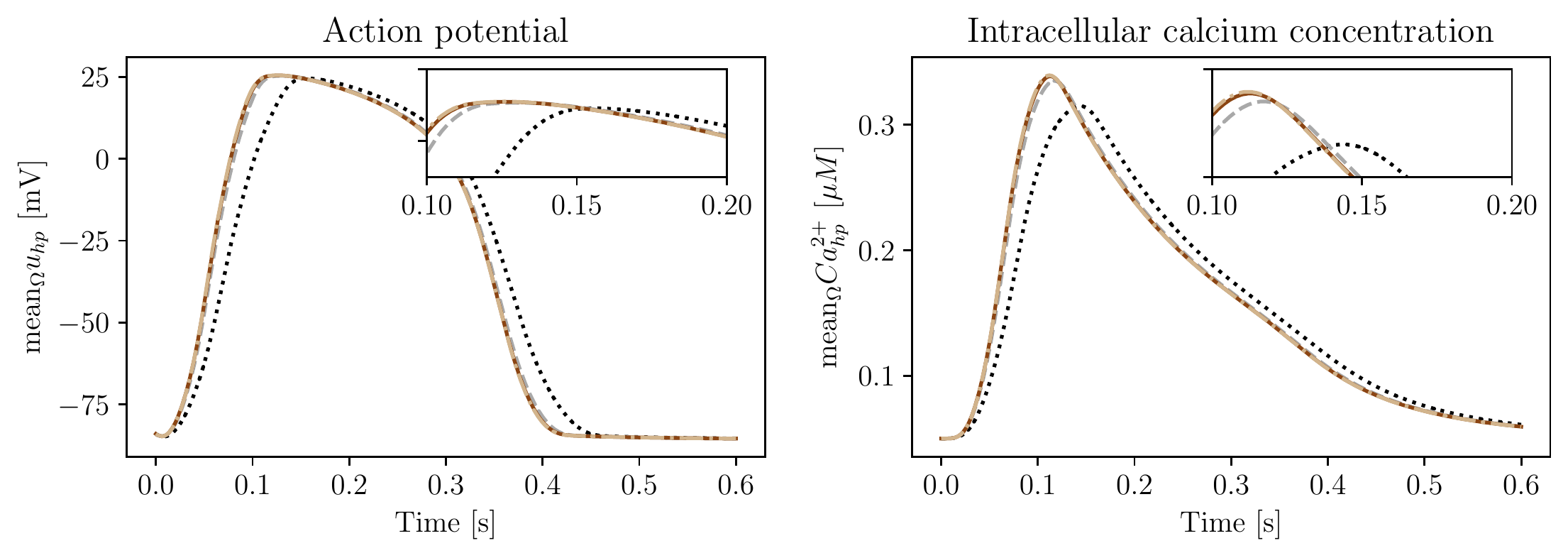}
    \end{subfigure}

    \hfill

    \begin{subfigure}{.95\linewidth}
        \centering
        \includegraphics[keepaspectratio, width=\textwidth]{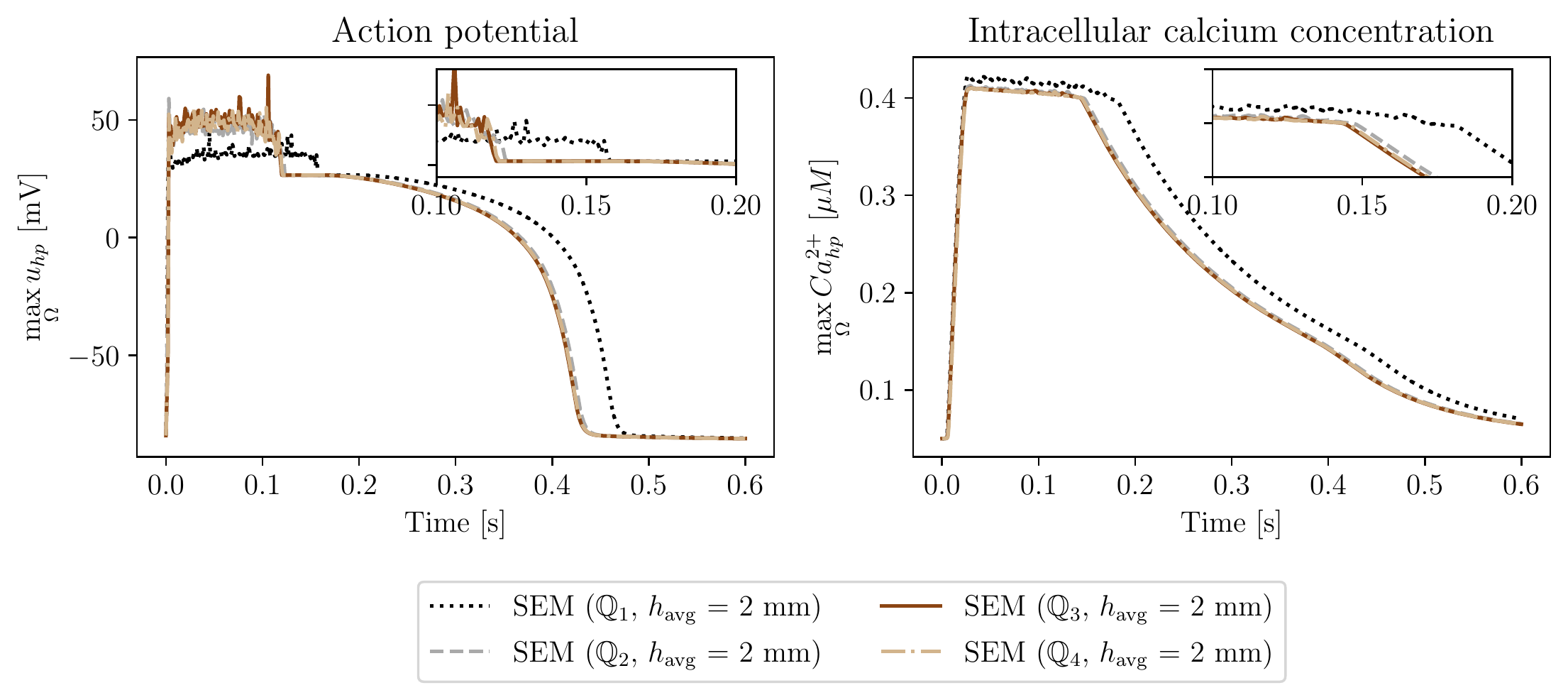}
    \end{subfigure}

    \caption{\emph{Zygote left ventricle}. Minimum (top), average (mid), maximum
        (bottom) pointwise values of action potential $u$ and intracellular calcium concentration $Ca^{2+}$ over time, ${\mathbb Q}_p$ with a fixed average mesh size $h_\mathrm{avg} = \SI{2.0}{\milli\meter}$ and $p$ ranging from 1 to 4.}
    \label{fig: LV_pointwise}
    \end{figure}

    In Figure~\ref{fig:activation_time_LV} we depict the activation maps for different choices of the
    local space (from ${\mathbb Q}_1$ to ${\mathbb Q}_4$). By looking at the
    contour lines, we observe that the ${\mathbb Q}_3$ solution is very close to the
    ${\mathbb Q}_4$ solution, that means we reach convergence for $p = 3$, even
    with such a relatively low mesh resolution $h_{\mathrm{avg}} =
    \SI{2.0}{\milli\meter}$. Whereas, it is a well--established result in the
    literature that first order Finite Elements would reach convergence for a value
    of $h_{\mathrm{avg}}$ that is about $100$ times smaller -- \textit{i.e.} for a
    much higher number of DOFs (see, \textit{e.g.}, \cite{Woodworth_2021}). We remark that here ``DOFs''
    refers to the number of degrees of freedom associated with the action potential, disregarding both
    gating variables and ionic species, thus it coincides with the number of mesh nodes. The same conclusions hold when considering Figure~\ref{fig: LV_pointwise},
    where we show the minimum, average, and maximum pointwise values of
    both the action potential $u$ and the intracellular calcium concentration $Ca^{2+}$ over time.

\begin{figure}
    \centering
    \captionsetup{justification=centering,margin=2cm}
    \scalebox{0.8}{
    \begin{subfigure}[t]{\linewidth}
        \centering
        \includegraphics[width=\textwidth]{legend_LV.png}
    \end{subfigure}
    }
    \\[\baselineskip]
    \scalebox{0.8}{
    \begin{subfigure}[t]{0.3\linewidth}
        \centering
        \includegraphics[width=\textwidth]{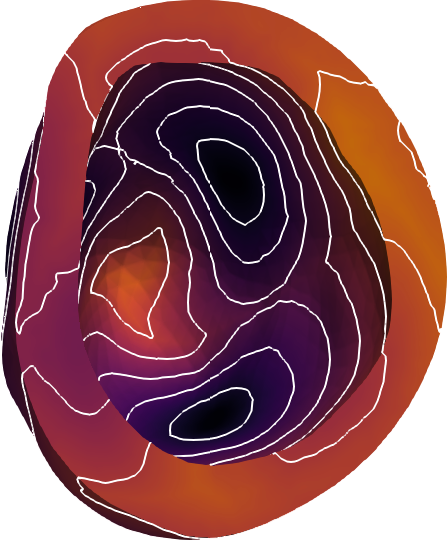}
        \caption*{}
    \end{subfigure}
    \hfill
    \begin{subfigure}[t]{0.3\linewidth}
        \centering
        \includegraphics[width=\textwidth]{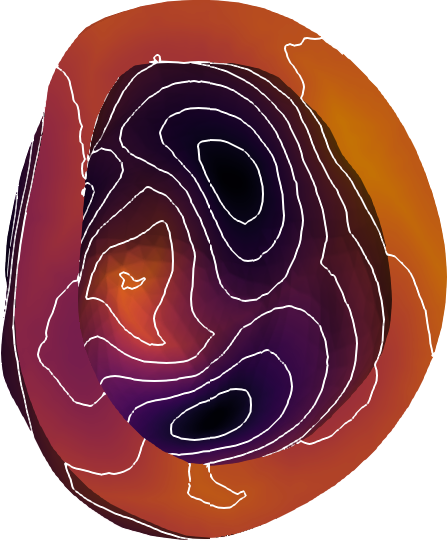}
        \caption*{}
    \end{subfigure}
    \hfill
    \begin{subfigure}[t]{0.3\linewidth}
        \centering
        \includegraphics[width=\textwidth]{output_view1_activation_time_EP_LV_Q4_0ref.png}
        \caption*{}
    \end{subfigure}
    }
    \\[\baselineskip]
    \scalebox{0.8}{
    \begin{subfigure}[t]{0.3\linewidth}
        \centering
        \includegraphics[width=\textwidth]{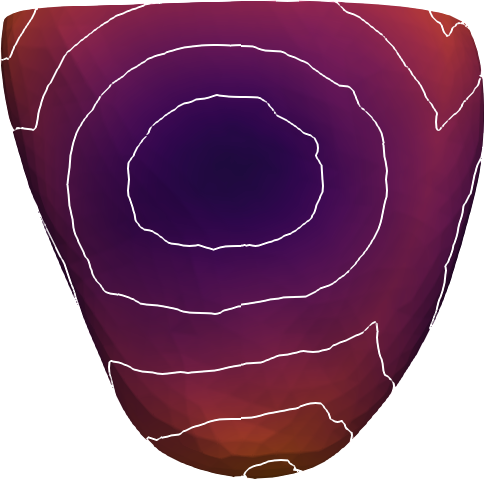}
        \caption*{}
    \end{subfigure}
    \hfill
    \begin{subfigure}[t]{0.3\linewidth}
        \centering
        \includegraphics[width=\textwidth]{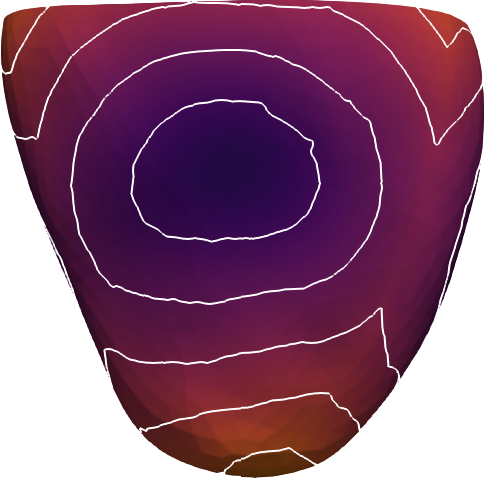}
        \caption*{}
    \end{subfigure}
    \hfill
    \begin{subfigure}[t]{0.3\linewidth}
        \centering
        \includegraphics[width=\textwidth]{output_view2_activation_time_EP_LV_Q4_0ref.png}
        \caption*{}
    \end{subfigure}
    }
    \\[\baselineskip]
    \scalebox{0.8}{
    \begin{subfigure}[t]{0.3\linewidth}
        \centering
        \includegraphics[width=\textwidth]{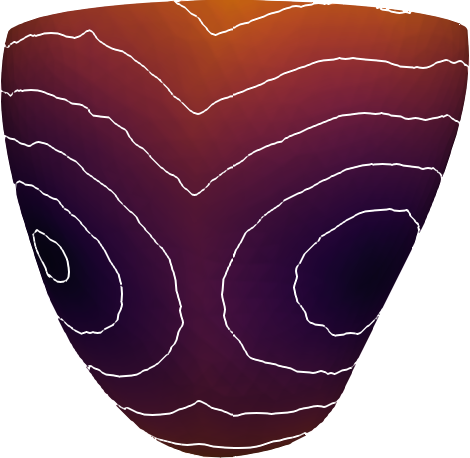}
        \caption*{$9'116'741 \,\, \text{DOFs}$\\($\mathbb{Q}_1$, $h_\mathrm{avg} = \SI{0.5}{\milli\meter}$)}
    \end{subfigure}
    \hfill
    \begin{subfigure}[t]{0.3\linewidth}
        \centering
        \includegraphics[width=\textwidth]{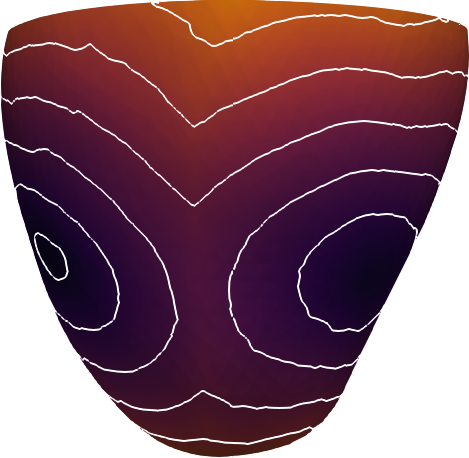}
        \caption*{$9'116'741 \,\, \text{DOFs}$\\($\mathbb{Q}_2$, $h_\mathrm{avg} = \SI{1.0}{\milli\meter}$)}
    \end{subfigure}
    \hfill
    \begin{subfigure}[t]{0.3\linewidth}
        \centering
        \includegraphics[width=\textwidth]{output_view3_activation_time_EP_LV_Q4_0ref.png}
        \caption*{$9'116'741 \,\, \text{DOFs}$\\($\mathbb{Q}_4$, $h_\mathrm{avg} = \SI{2.0}{\milli\meter}$)}
    \end{subfigure}
    }
    \captionsetup{justification=raggedright,margin=0cm}
    \caption{\emph{Zygote left ventricle}. Three views of the activation maps computed with ${\mathbb Q}_1, {\mathbb Q}_2, {\mathbb Q}_4$ elements and a fixed number of DOFs equal to $9'116'741$.}
    \label{fig:activation_time_LV_dofs}
    \end{figure}

\begin{figure}
    \centering
    \begin{subfigure}{.95\linewidth}
        \centering
        \includegraphics[keepaspectratio, width=\textwidth]{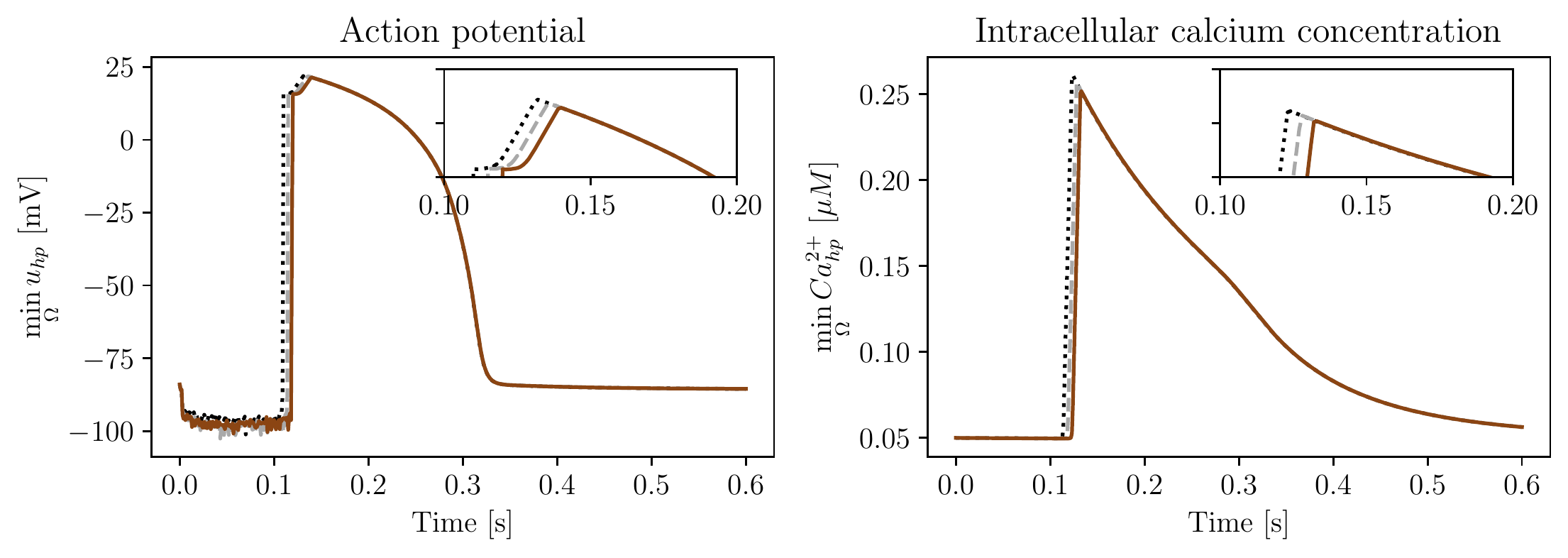}
    \end{subfigure} %

    \hfill

    \begin{subfigure}{.95\linewidth}
        \centering
        \includegraphics[keepaspectratio, width=\textwidth]{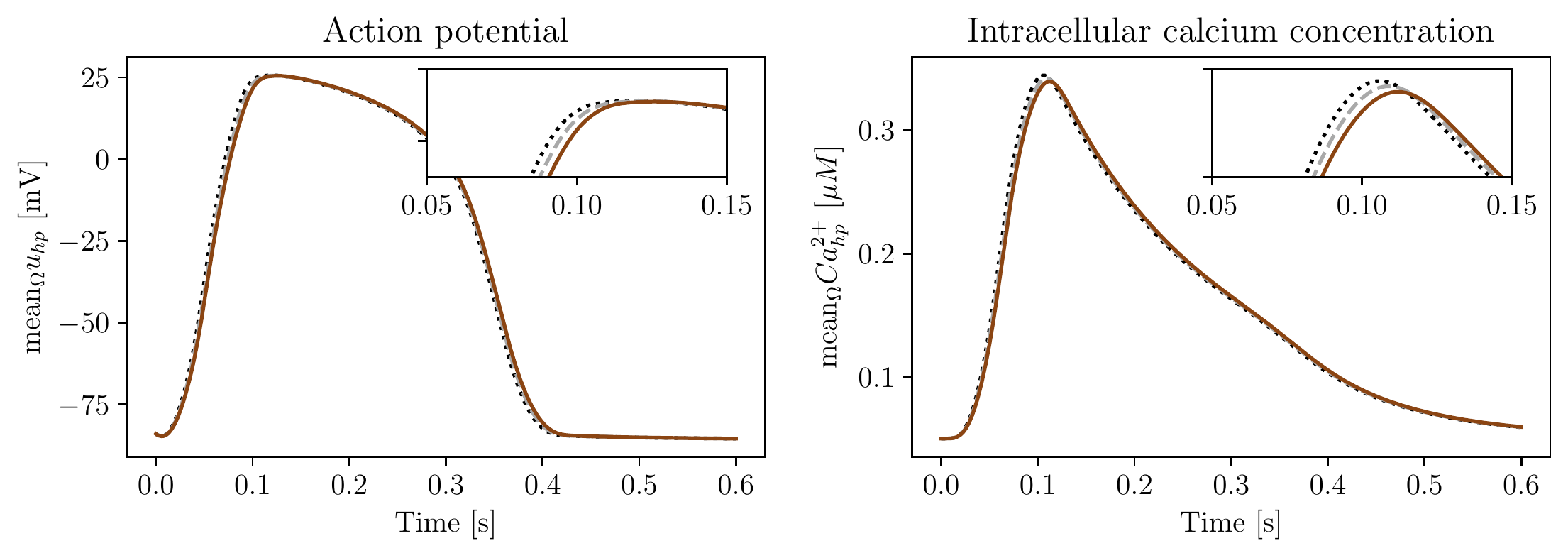}
    \end{subfigure}

    \hfill

    \begin{subfigure}{.95\linewidth}
        \centering
            \includegraphics[keepaspectratio, width=\textwidth]{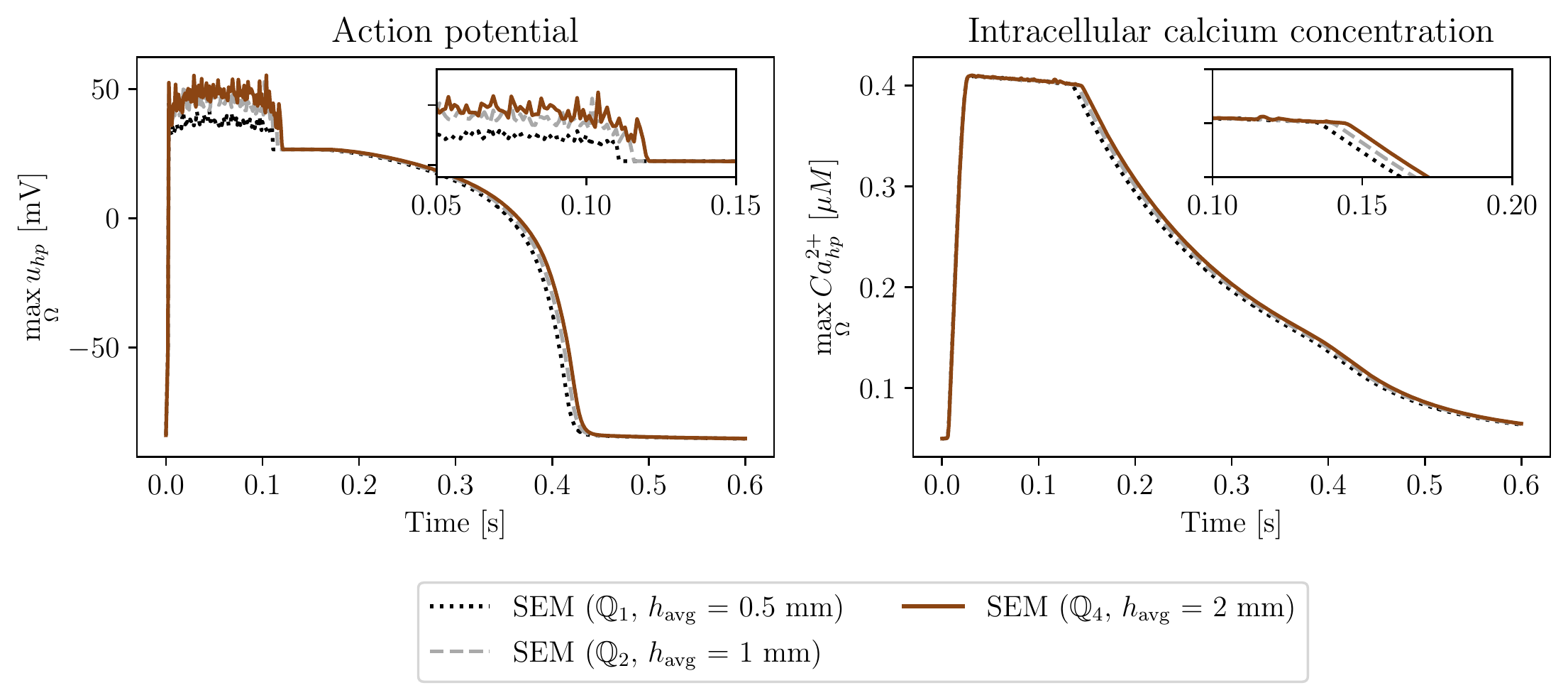}
    \end{subfigure}

    \caption{\emph{Zygote left ventricle}. Minimum (top), average (mid), maximum
        (bottom) pointwise values of action potential $u$ and intracellular calcium concentration $Ca^{2+}$ over time, for ${\mathbb Q}_1, {\mathbb Q}_2, {\mathbb Q}_4$ elements and a fixed number of DOFs equal to $9'116'741$.}
    \label{fig: LV_pointwise_dofs}
    \end{figure}

    To further verify these conclusions, we consider different combinations of \(h\) and \(p\) that lead to the same number of DOFs, namely \(9'116'741\) for $\mathbb{Q}_1$ with $h_\mathrm{avg} = \SI{0.5}{\milli\meter}$, $\mathbb{Q}_2$ with $h_\mathrm{avg} = \SI{1.0}{\milli\meter}$ and $\mathbb{Q}_4$ with $h_\mathrm{avg} = \SI{2.0}{\milli\meter}$. The activation maps displayed in Figure~\ref{fig:activation_time_LV_dofs} and the pointwise values of the action potential and intracellular calcium shown in Figure~\ref{fig: LV_pointwise_dofs} reveal that all the results are quite similar and pretty close to convergence, with the \(\mathbb{Q}_4\) simulation still being the most accurate. Table~\ref{tab:LV_timing_dofs} summarizes the parameters and the computational times recorded for the three simulations, which have been performed with the matrix--free solver and the SEM formulation.
    Going from $\mathbb{Q}_1$ to $\mathbb{Q}_4$, the time spent in solving the linear system is reduced of about 9\%, whereas the cost of assembling the right–hand side is reduced of about 50\%, leading to an overall reduction of about 12\%. These results further confirm that in the matrix--free context the strategy of increasing \(p\) rather than reducing \(h\) is more advantageous in terms of both numerical accuracy and computational efficiency.

    \begin{table}
    \hspace*{-.8cm}
    {\fontsize{8}{9}\selectfont{
            \begin{tabular}{ | c c c c | c c c |}
                \toprule
                \makecell{Mesh points \\ number }  & \makecell{Cells \\ number} & \makecell{$h_\mathrm{avg}$ \\ \text{[\si{\milli\meter}]}} & \makecell{Local \\ space} & \makecell{PCG iterations\\ GMG preconditioner} & \makecell{Linear solver \\ {} [\si{\second}]} & \makecell{Assemble rhs \\ {}[\si{\second}]} \\
                \midrule
                $9'116'741$ & $8'939'776$ & 0.05 & $\mathbb{Q}_1$ & 2.54 & 12440.968 & 1241.884 \\
                $9'116'741$ & $1'117'472$ & 0.1  & $\mathbb{Q}_2$ & 2.54 & 12334.152 & 1020.451 \\
                $9'116'741$ &   $139'684$ & 0.2  & $\mathbb{Q}_4$ & 2.07 & 11328.664 &  687.311 \\
                \bottomrule
            \end{tabular}
    }}
    \caption{\emph{Zygote left ventricle}. Computational times for ${\mathbb Q}_1, {\mathbb Q}_2, {\mathbb Q}_4$ elements and a fixed number of DOFs equal to $9'116'741$.}
    \label{tab:LV_timing_dofs}
    \end{table}

    \subsection{Whole--heart}

    The aim of this section is to show that our matrix--free solver can be
    successfully applied even in a much more complex framework. For this purpose we perform
    a numerical simulation in sinus rhythm with the Zygote four--chamber heart \cite{Zygote}.
    The settings of this test case can be found at the beginning of Section~\ref{sec:numres}.
    We consider different ionic models, namely CRN \cite{courtemanche1998} and TTP06 \cite{tentusscher_2006} for atria and ventricles, respectively.
    Furthermore, we model the valvular rings as non-conductive regions of the myocardium.
    The mesh is endowed with 355'664 cells and 10'355'058 nodes ($h_{\mathrm{avg}} = \SI{1.6}{\milli\meter}$).
    We employ the matrix--free solver and ${\mathbb Q}_3-$SEM,
    this choice is motivated by the numerical results obtained for the Niederer benchmark (Section~\ref{sec:slab}) and the convergence test performed on the Zygote left ventricle geometry (Section~\ref{sec:lv}).

\begin{figure}
        \centering
        \captionsetup{justification=centering,margin=2cm}
        \begin{subfigure}[t]{\linewidth}
            \centering
            \includegraphics[width=\textwidth]{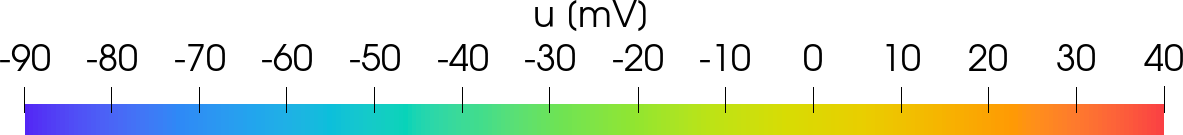}
        \end{subfigure}
        \\[\baselineskip]
        \begin{subfigure}[t]{0.2\linewidth}
            \centering
            \includegraphics[width=\textwidth]{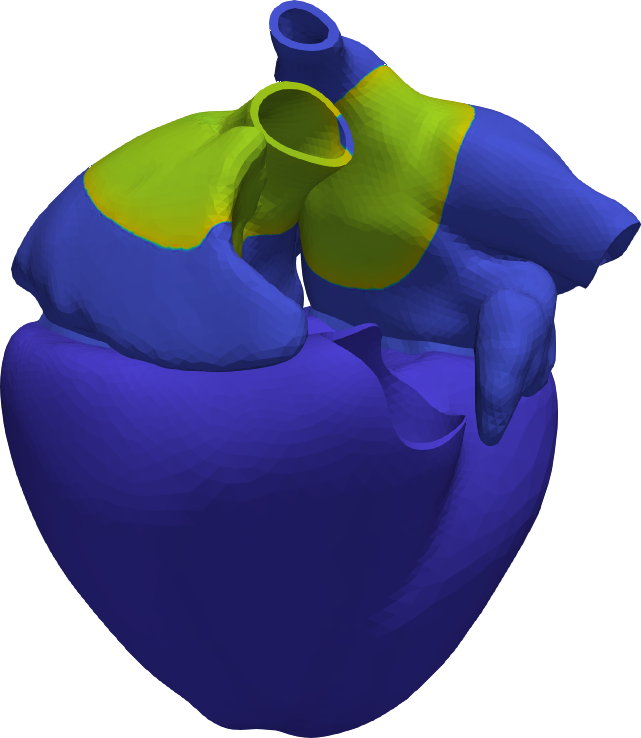}
            \caption*{$t = \SI{0.1}{\second}$}
        \end{subfigure}
        \hfill
        \begin{subfigure}[t]{0.2\linewidth}
            \centering
            \includegraphics[width=\textwidth]{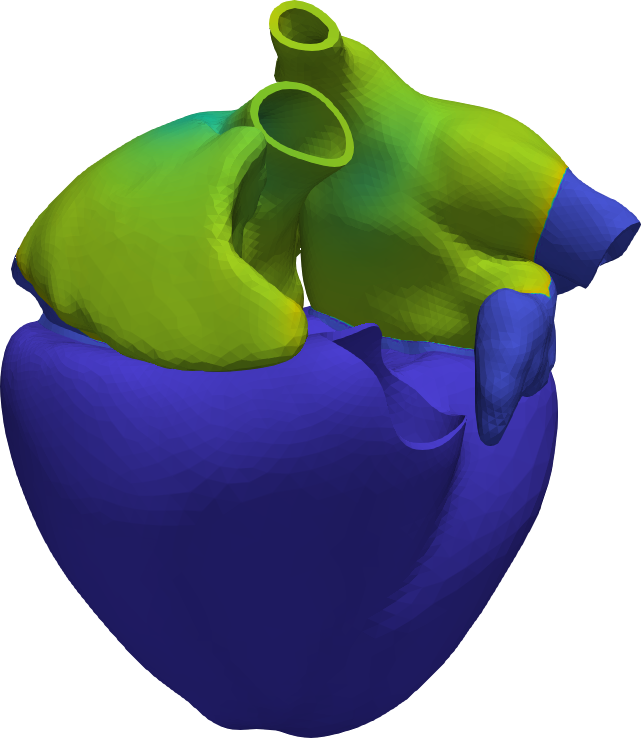}
            \caption*{$t = \SI{0.2}{\second}$}
        \end{subfigure}
        \hfill
        \begin{subfigure}[t]{0.2\linewidth}
            \centering
            \includegraphics[width=\textwidth]{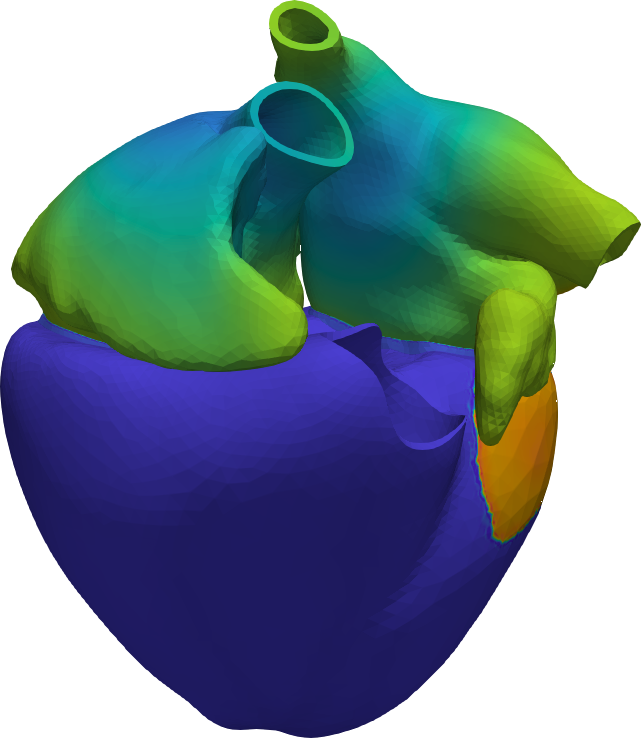}
            \caption*{$t = \SI{0.3}{\second}$}
        \end{subfigure}
        \hfill
        \begin{subfigure}[t]{0.2\linewidth}
            \centering
            \includegraphics[width=\textwidth]{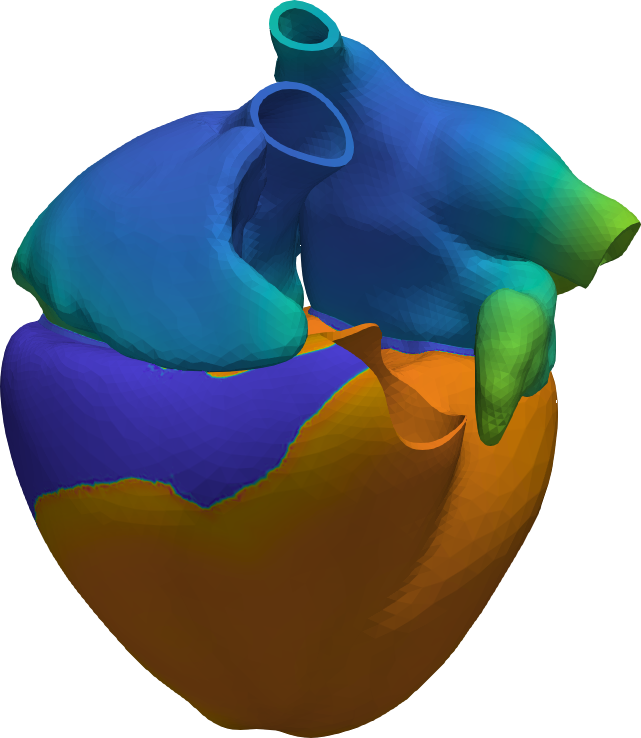}
            \caption*{$t = \SI{0.4}{\second}$}
        \end{subfigure}
        \\
        \begin{subfigure}[t]{0.2\linewidth}
            \centering
            \includegraphics[width=\textwidth]{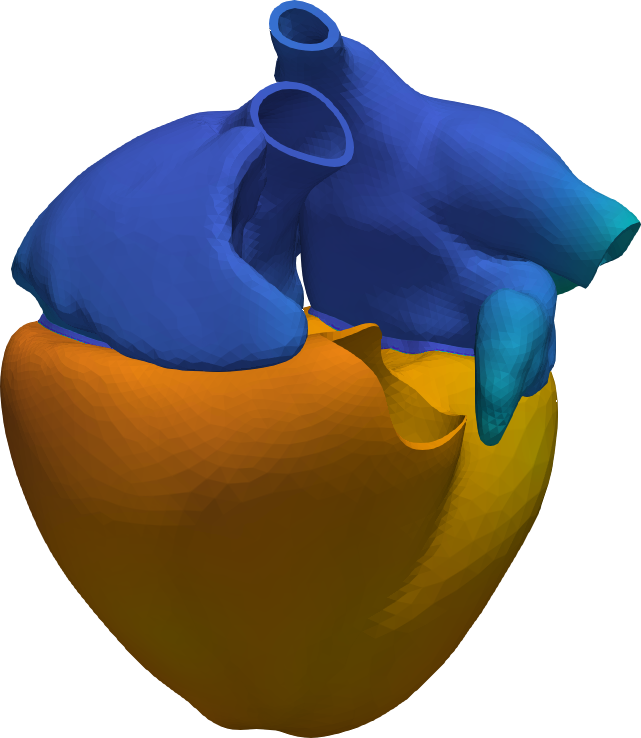}
            \caption*{$t = \SI{0.5}{\second}$}
        \end{subfigure}
        \hfill
        \begin{subfigure}[t]{0.2\linewidth}
            \centering
            \includegraphics[width=\textwidth]{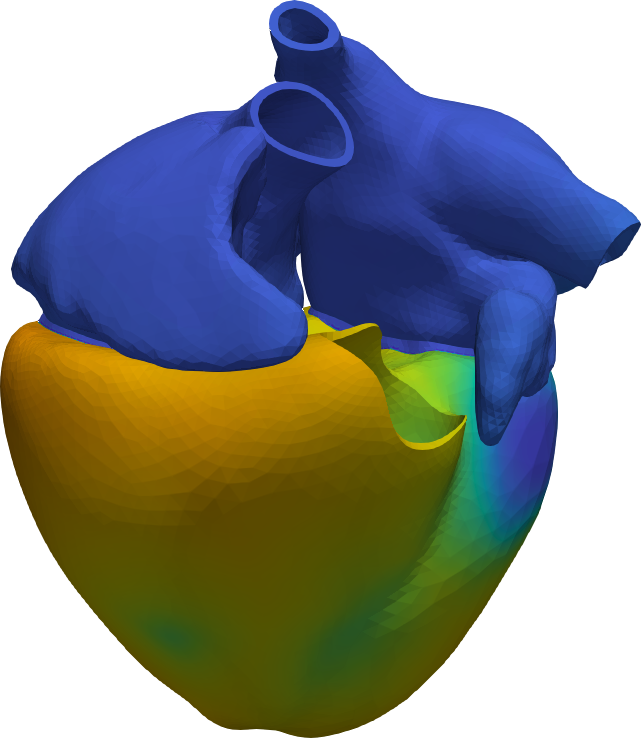}
            \caption*{$t = \SI{0.6}{\second}$}
        \end{subfigure}
        \hfill
        \begin{subfigure}[t]{0.2\linewidth}
            \centering
            \includegraphics[width=\textwidth]{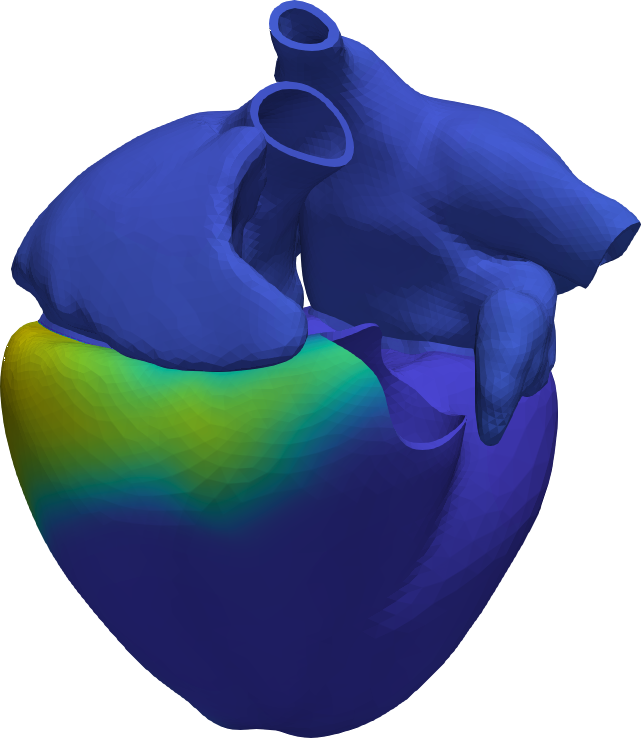}
            \caption*{$t = \SI{0.7}{\second}$}
        \end{subfigure}
        \hfill
        \begin{subfigure}[t]{0.2\linewidth}
            \centering
            \includegraphics[width=\textwidth]{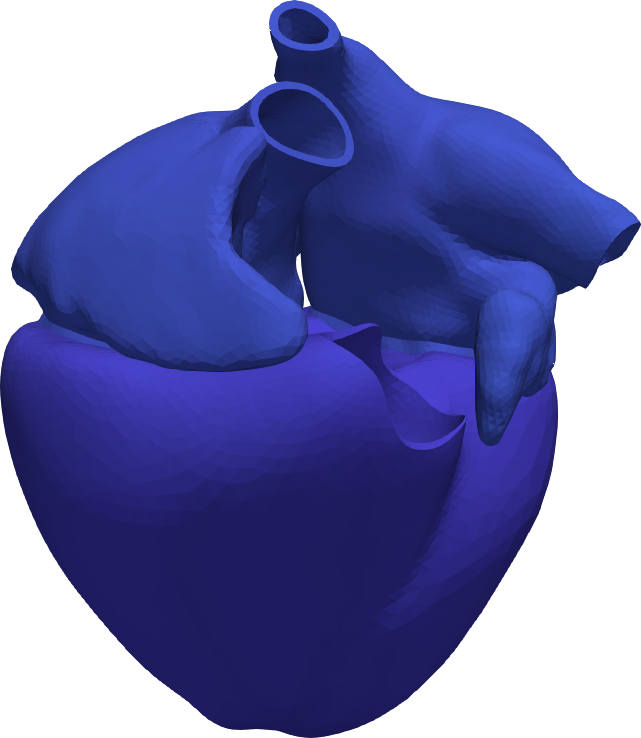}
            \caption*{$t = \SI{0.8}{\second}$}
        \end{subfigure}
        \captionsetup{justification=raggedright,margin=0cm}
        \caption{\emph{Zygote whole--heart}. Time evolution of the transmembrane potential $u$ during one heartbeat.}
        \label{fig:solution_4CH}
    \end{figure}

    We depict in Figure~\ref{fig:solution_4CH} the evolution of the transmembrane
    potential over time on the whole--heart geometry.
    The electric signal initiates at the sinoatrial node in the right atrium and then
    propagates to the left atrium and ventricles by means of preferential conduction lines,
    such as the Bachmann's and His bundles \cite{Harrington_2011}.
    The wavefront propagation appears very smooth, while accounting for small dissipation and dispersion throughout the heartbeat, as expected from the use of high--order discretizations \cite{Bucelli_2021}.

    \section{Conclusions}
    We developed a matrix--free solver for cardiac electrophysiology
    tailored to the efficient use of high--order numerical methods.
    We employed the monodomain equation to model the propagation of
    the transmembrane potential and physiologically--based ionic models
    (CRN and TTP06) to describe the behaviour of different chemical species at the
    cell level.

    We run several electrophysiological simulations for three different test cases, namely a slab of cardiac tissue,
    the Zygote left ventricle and the Zygote whole--heart to demonstrate the effectiveness
    and generality of our matrix--free solver in combination with Spectral Element Methods.
    SEM and SEM--NI provided comparable numerical results in terms of both accuracy and efficiency.
    Furthermore, we showed the importance of considering high--order Finite Elements in improving the accuracy and the computational burden for this class
    of mathematical problems, \textit{i.e.} with sharp wavefronts involved.

    Our matrix--free solver outperforms state--of--the--art matrix--based solvers
    in terms of computational costs and memory requirements.
    This is true even when matrix--vector products are computed without any matrix--free preconditioner,
    thanks to both vectorization and sum--factorization.
    Finally, the low memory footprint of the matrix--free implementation may allow for the development of
    GPU--based solvers of the cardiac function.

\appendix
\section{}\label{Appendix}
Let us consider
the domain $\Omega\times(0,T)$ with $\Omega=(0,\SI{1}{\centi\meter})^2$ and
$T=\SI{1}{\second}$. For the discretization we have taken a mesh size $h=\SI{0.5}{\milli\meter}$
(consistently with those chosen in the simulations shown in Section~\ref{sec:numres}),
polynomial degree $p=3$ and $p=4$, and the BDF schemes of order 1, 2, and 3 in
time. Then we have measured the errors

\begin{subequations}\label{eq:errors-for-BDF}
\begin{align}
\mathrm{err}_{H^1} & =\left(\Delta t \sum_n \left\|u_{hp}(t_n)-u(t^n)\right\|^2_{H^1(\Omega)}\right)^{\frac{1}{2}}, \\
\mathrm{err}_{L^2} & =\left(\Delta t \sum_n \left\|u_{hp}(t_n)-u(t^n)\right\|^2_{L^2(\Omega)}\right)^{\frac{1}{2}},
\end{align}
\end{subequations}
between the exact solution $u$ and the fully discrete SEM--NI solution $u_{hp}$.

\begin{figure}
        \centering
        \includegraphics[keepaspectratio, width=\textwidth]{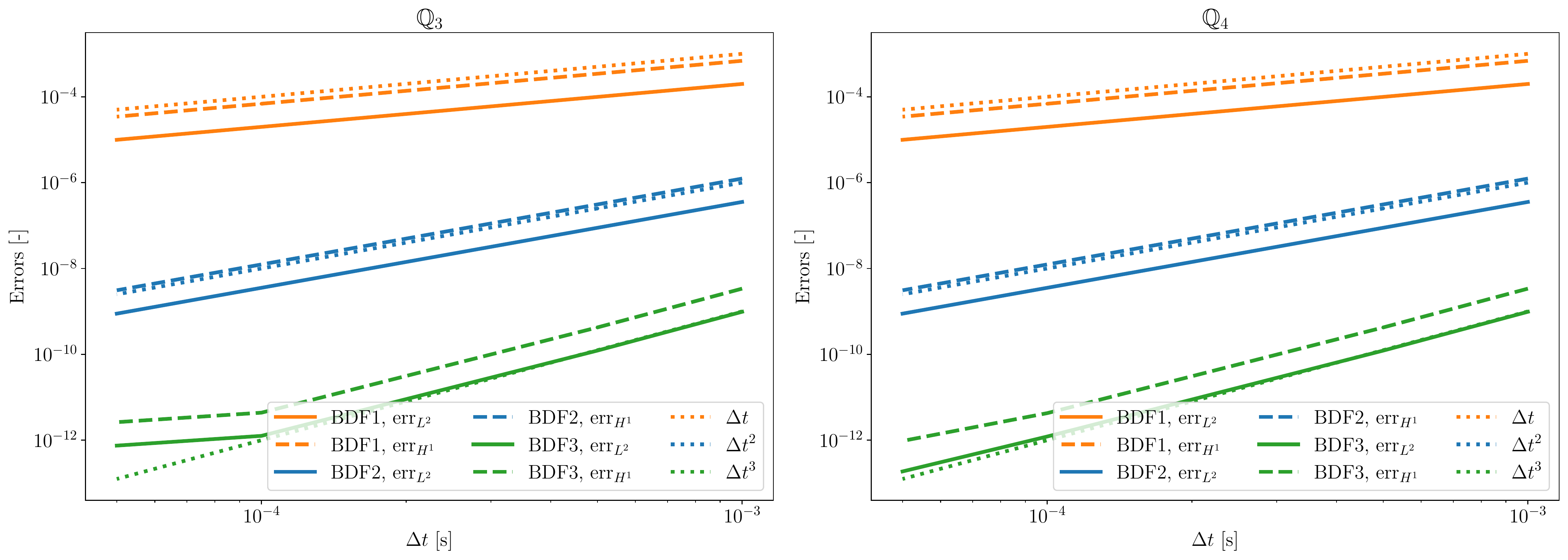}
        \caption{Errors \eqref{eq:errors-for-BDF} for SEM--NI
discretization of the heat equation, with $h=0.5$mm and test solution
$u({\bf x},t)=(x_1^2+x_2^2)(2+\sin(\pi t))$. $p=3$ on the left and $p=4$ on the
right}
        \label{fig:bdf-accuracy1}
    \end{figure}

\begin{figure}
        \centering
        \includegraphics[keepaspectratio, width=\textwidth]{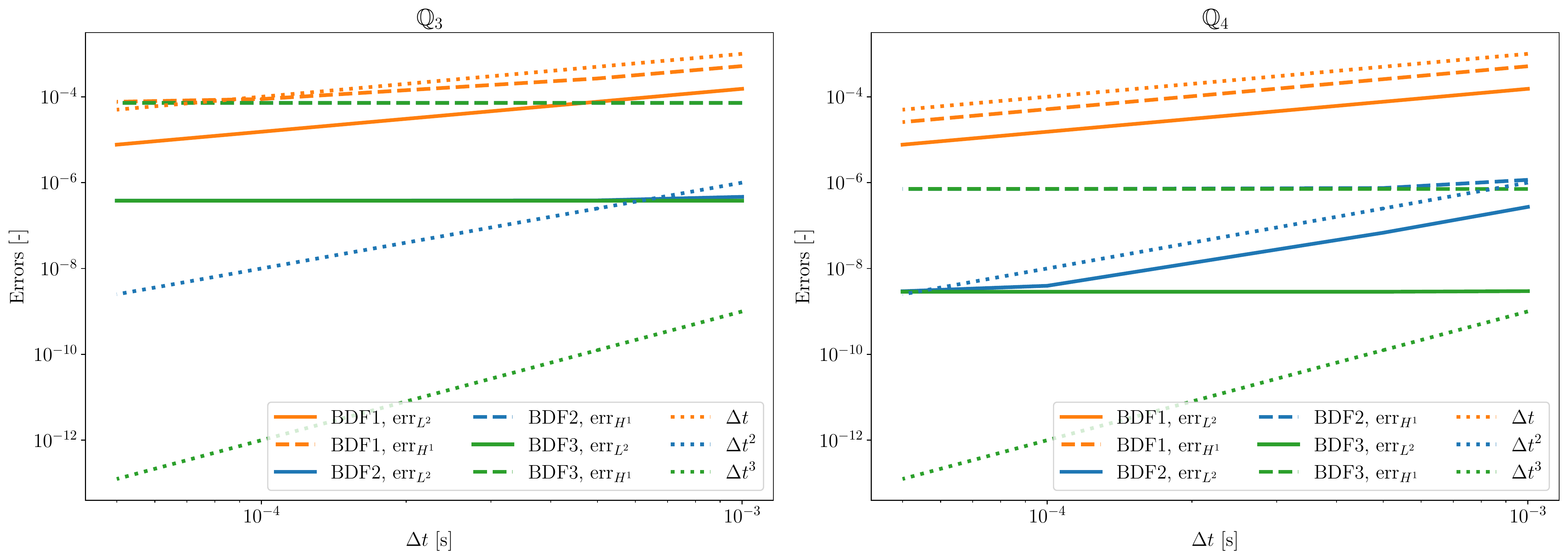}
        \caption{Errors (\ref{eq:errors-for-BDF}) for SEM--NI
discretization of the heat equation, with $h=0.5$mm and test solution $u({\bf
x},t)=\sin(\pi x_1)\sin(\pi x_2)(2+\sin(\pi t))$.  $p=3$ on the left and $p=4$ on the
right}
        \label{fig:bdf-accuracy2}
    \end{figure}

We have considered the exact solutions
$u({\bf x},t)=(x_1^2+x_2^2)(2+\sin(\pi t))$
and $u({\bf x},t)=\sin(\pi x_1)\sin(\pi x_2)(2+\sin(\pi t))$.
The former one is solved exactly in space by both ${\mathbb Q}_3$ and ${\mathbb
Q}_4$ discretizations, so that in Figure~\ref{fig:bdf-accuracy1} we can appreciate the full convergence order in
time of all the three schemes BDF1, BDF2, and BDF3.
On the contrary, the latter solution is
not captured exactly. The associated errors are shown in  Figure~\ref{fig:bdf-accuracy2}
and we observe that they are bounded from
below from the space discretization error when $\Delta t\lesssim 10^{-4}\,\si{\second}$,
making the higher accuracy of BDF3 worthless.
We notice that $\Delta t=10^{-4}\,\si{\second}$ coincides with the time step of $\SI{0.1}{\milli\second}$ used in
the simulations reported in Section~\ref{sec:numres}.
We remark that the numerical solutions we are looking for in cardiac electrophysiology
typically feature steepest gradients than those we have considered
in these examples, thus it is  unlikely that the plateaux in the errors starts in
correspondence of time steps much smaller than $10^{-4}\,\si{\milli\second}$.
Finally, we bear in mind that, while BDF2 is absolutely stable, BDF3 is
not, then special care should be given to the choice of the time step.

    \section*{Acknowledgments}
    This research has been funded partly by the Italian Ministry of University and Research
    (MIUR) within the PRIN (Research projects of relevant national interest 2017
    ``Modeling the heart across the scales: from cardiac cells to the whole organ''
    Grant Registration number 2017AXL54F).

    \bibliography{bibliography}
\end{document}

%% file: packages.tex
\usepackage[english]{babel}
\usepackage[utf8]{inputenc}

\usepackage[hyperref]{xcolor}
\definecolor{lifex}{HTML}{f60248}
\newcommand{\lifex}{\texorpdfstring{\texttt{life\textsuperscript{\color{lifex}{x}}}}{lifex}}

\usepackage[hidelinks]{hyperref}

\usepackage{amsmath}
\usepackage{amssymb}
\usepackage{amsthm}
\usepackage[overload]{empheq}
\usepackage{bm} 

\DeclareMathOperator*{\mean}{mean}

\usepackage{siunitx}
\usepackage{graphicx}
\graphicspath{{figures/}}
\usepackage{subcaption}

\usepackage{tabularx}
\usepackage{makecell}
\usepackage{booktabs}

\usepackage[linesnumbered]{algorithm2e}
\SetEndCharOfAlgoLine{} 
\SetArgSty{upshape} 

\let\oldnl\nl
\newcommand{\nonl}{\renewcommand{\nl}{\let\nl\oldnl}}

\DeclareMathOperator*{\argmax}{arg\,max}

\usepackage{etoolbox}
\makeatletter
\patchcmd\ps@pprintTitle
{\fi\fi\fi\fi}
{\fi\fi\fi\fi\afterassignment\fix@elsarticle\let\@tempa}{}{\FailedToPatch}
\def\fix@elsarticle{\iffalse{\fi}\romannumeral-`0}
\makeatother

%% file: macros.tex
\newcommand{\papertitle}{A matrix--free high--order solver for the numerical
solution of cardiac electrophysiology}

\newcommand{\fZero}{{\mathbf{f}_0}}
\newcommand{\sZero}{{\mathbf{s}_0}}
\newcommand{\nZero}{{\mathbf{n}_0}}

\newcommand{\Chi}{\chi}
\newcommand{\Cm}{C_\mathrm{m}}
\newcommand{\Iion}{{\mathcal{I}_{\mathrm{ion}}}}
\newcommand{\Iapp}{{\mathcal{I}_{\mathrm{app}}}}
\newcommand{\DiffTens}{\boldsymbol{D}_{\mathrm{M}}}
\newcommand{\Pot}{u}
\newcommand{\Gating}{\boldsymbol{w}}
\newcommand{\Conc}{\boldsymbol{z}}
\newcommand{\RhsGating}{\boldsymbol{H}}
\newcommand{\RhsConc}{\boldsymbol{G}}
\newcommand{\sigmal}{\sigma_{\text{l}}}
\newcommand{\sigmat}{\sigma_{\text{t}}}
\newcommand{\sigman}{\sigma_{\text{n}}}

\newcommand{\APot}{\mathsf u}
\newcommand{\AGating}{\bm{\mathsf w}}
\newcommand{\AConc}{\bm{\mathsf z}}
\newcommand{\AIion}{\mathsf s}
\newcommand{\AIapp}{\mathsf f}

\newcommand{\Mass}{\mathsf M}
\newcommand{\Stiff}{\mathsf K}
\newcommand{\Matrix}{\mathsf A}
\newcommand{\ARhs}{\mathsf b}